\documentclass[11pt]{amsart}

\usepackage{graphicx}
\usepackage{mathptmx}      

\usepackage[english]{babel}
\usepackage[latin1]{inputenc} 
\usepackage{epstopdf}
\usepackage{verbatim} 
\usepackage{amsmath,amssymb,amsfonts,mathrsfs,latexsym} 
\usepackage{fancyhdr}
\usepackage{xr}

\newcommand{\R}{{\mathbb R}}

\newcommand{\dd}{\mathrm{d}}

\newcommand{\diff}[2]{\frac{\dd #1}{\dd #2}}

\DeclareMathOperator{\sn}{sn}
\DeclareMathOperator{\cn}{cn}
\DeclareMathOperator{\dn}{dn}
\DeclareMathOperator{\am}{am}

\newcommand{\minimize}{\operatornamewithlimits{minimize}}
\newcommand{\bE}{\mathcal{E}}
\newcommand{\fr}{\colon}
\renewcommand{\a}{\alpha}
\renewcommand{\b}{\beta}
\newcommand{\y}{\boldsymbol{\gamma}}
\renewcommand{\l}{\lambda}
\renewcommand{\k}{\kappa}
\newcommand{\w}{S}
\renewcommand{\th}{\theta}
\newcommand{\bbn}{\mathbf{n}}
\newcommand{\bbt}{\mathbf{t}}

\newcommand{\infl}{\boldsymbol{\zeta}}

\newcommand{\stp}{s_0}
\newcommand{\dtheta}{\mathrm{d}\theta}
\newcommand{\ds}{\mathrm{d}s}
\newcommand{\dt}{\mathrm{d}t}
\newcommand{\du}{\mathrm{d}u}

\newcommand{\bbx}{\mathbf{x}}

   \newtheorem{theorem}{Theorem}[section]

 \theoremstyle{remark}

   \newtheorem*{note}{Remark}
   
\numberwithin{equation}{section}

\begin{document}

\title{Approximation by planar elastic curves
}
 \thanks{Research supported by Innovation Fund Denmark, project number  128-2012-3}

\author{David Brander}
\address{Department of Applied Mathematics and Computer Science,\\
Technical University of Denmark\\
 Matematiktorvet, Building 303 B\\
DK-2800 Kgs. Lyngby\\ Denmark}
\email{dbra@dtu.dk}

\author{Jens Gravesen}
\address{Department of Applied Mathematics and Computer Science,\\
Technical University of Denmark\\
 Matematiktorvet, Building 303 B\\
DK-2800 Kgs. Lyngby\\ Denmark}
\email{jgra@dtu.dk}

\author{Toke Bjerge N\o{}rbjerg}
\address{Department of Applied Mathematics and Computer Science,\\
Technical University of Denmark\\
 Matematiktorvet, Building 303 B\\
DK-2800 Kgs. Lyngby\\ Denmark}
\email{tono@dtu.dk}


\begin{abstract}
We give an algorithm for approximating a given plane curve segment by a planar elastic curve.
The method depends on an analytic representation of the space of elastic curve segments,
together with a geometric method for obtaining a good initial guess for the approximating
curve.  A gradient-driven optimization is then used to find the approximating elastic curve.
\end{abstract}

\maketitle

\section{Introduction}
 An \emph{Euler elastica} or \emph{elastic curve} is the solution to the
variational problem of minimizing the bending energy,
the integral of the curvature squared $\int \kappa(s)^2 \ds$,
among
curves of a given length with fixed endpoints and with the tangents prescribed at the
endpoints. All solutions to this problem
were described by Euler \cite{euler} in 1744, and the curves can be parameterized in terms of elliptic functions.

The energy minimizing property qualifies the elastica as a mathematical
model for the shape assumed by a thin inextensible rod when constraints
are place only at the endpoints. This shape appears
naturally in certain manufacturing scenarios: for example in a construction
made from thin, flexible strips of wood or similar, the shape
of each strip, between the fixed points, is an elastica.
In another application,  a thin metal blade can be heated and used to
cut polystyrene for architectural formwork,  the blade changing its shape during the
motion. This permits the construction of curved geometries potentially much more cheaply than
with the alternative of numerically controlled milling;  this article is part of a larger project dedicated
to the development of this so-called \emph{robotic hot-blade cutting} technology \cite{madifa}.
Here, too, the shape of the blade is an elastica.

 The present standard representation for curves and surfaces
 in computer-aided design (CAD) systems is \emph{rational} splines, that is, piecewise rational functions.
In architecture there is typically a second step after the initial conceptual design, whereby
the CAD model is adapted slightly for practical realization with respect to the chosen manufacturing
process, and this is called \emph{rationalization}.
 The problem of rationalizing a standard CAD design
for a construction method involving elastic curves entails the approximation
of a rational or polynomial spline curve by an elastica, and it is this problem
that motivates us here.   It is  worth noting that
the term ``spline'', in pre-CAD years, referred to a thin strip of wood used to draw
curves that interpolate a given set of points smoothly.  Between the interpolation
points, these strips assumed the shapes of elastica, or rather \emph{planar}
elastica, because the strips were laid on a flat surface. Thus, for the type of
rationalization referred to above, one could say that the task is, somewhat ironically, to approximate a
digital  spline with an analogue spline.

With this in mind,  we consider here the problem of approximating
an arbitrary curve by a planar elastica.\\

\noindent \textbf{Precursors:}
 Quite distinctly from the applications mentioned above, elastic curves are a good
candidate for the choice of functions used in geometric modeling because their curvature minimizing
property makes them optimally faired.   The idea of using them as mathematical splines -- here meaning piecewise
smooth functions -- has been discussed by many authors: for example the references
\cite{birkhoffdeboor1965}, \cite{borbelyjohnson2014}, \cite{golombJerome1982}, \cite{horn1983}, \cite{malcolm1977},
\cite{mehlum1974} and \cite{mumford1990}
constitute a representative, though not an exhaustive, list.
These works are concerned  either
with the problem of how to compute the elastic curve satisfying various constraints respecting placement of
endpoints, end-tangents, lengths etc., or with such problems as the existence of an elastic curve interpolating a given set of points.
As a testament to the importance of these curves, one finds that many aspects
of the  theory are re-derived independently several times.

For our purpose, however, the main point is that previous work on elastica and their
applications in spline theory does not
consider the problem of approximating an arbitrary curve by an elastica. \\

\begin{figure}[h]
\centering
\unitlength=0.25\textwidth
\includegraphics[width=\unitlength,trim=20mm 20mm 20mm 20mm,clip]{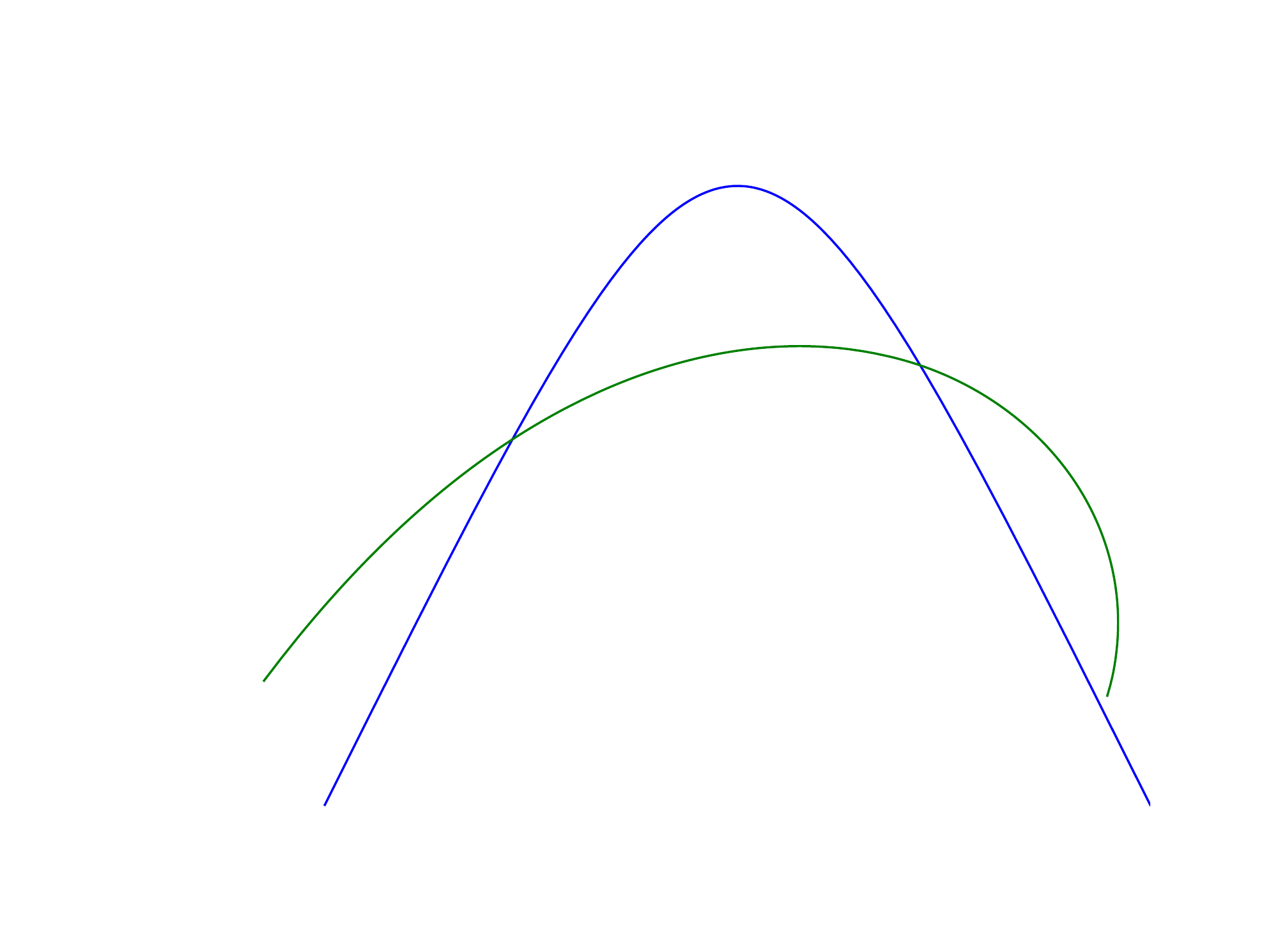}\hfil
\includegraphics[width=\unitlength,trim=35mm 30mm 35mm 30mm,clip]{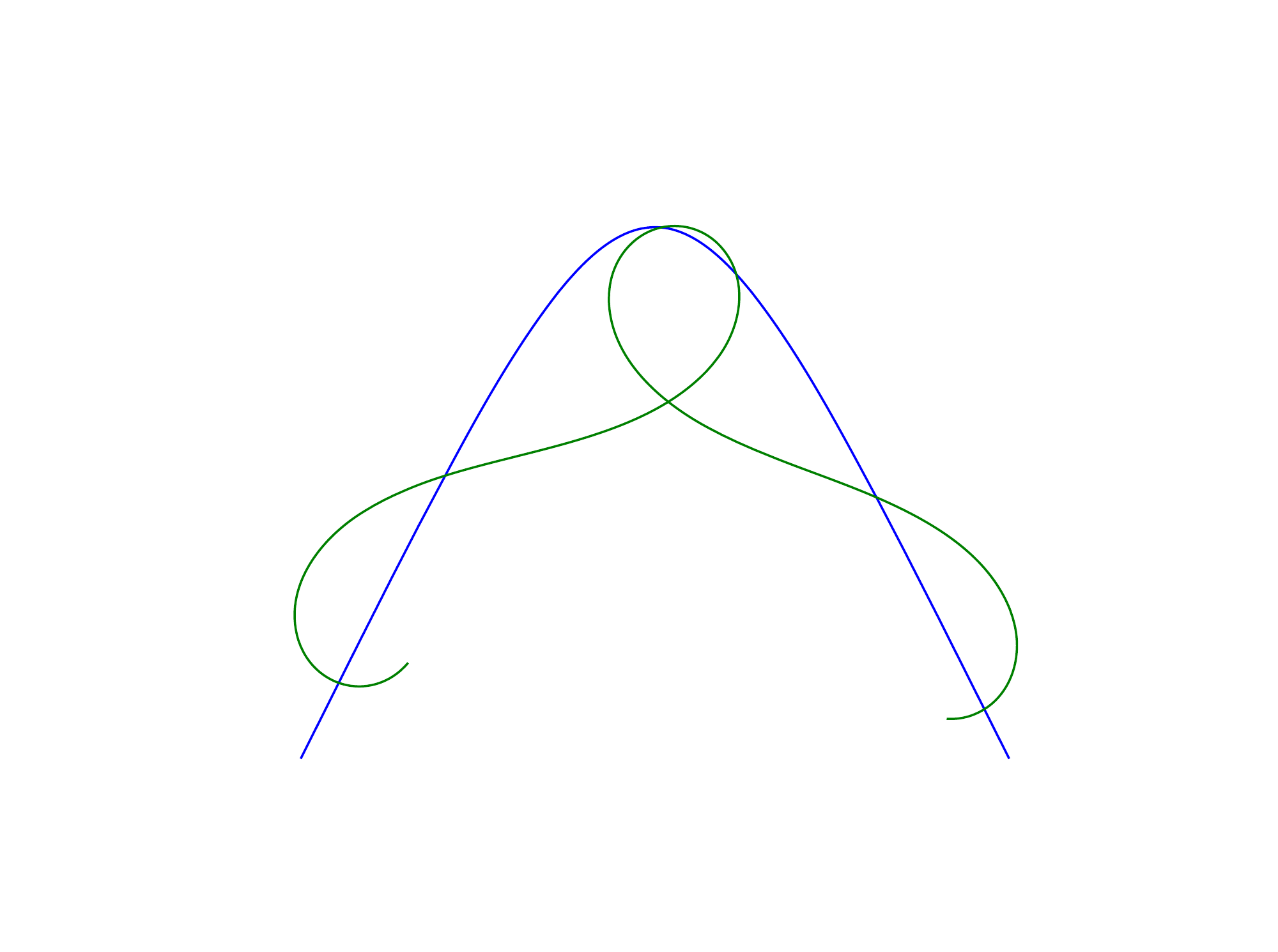}\hfil
\includegraphics[width=\unitlength,trim=19mm 20mm 19mm 20mm,clip]{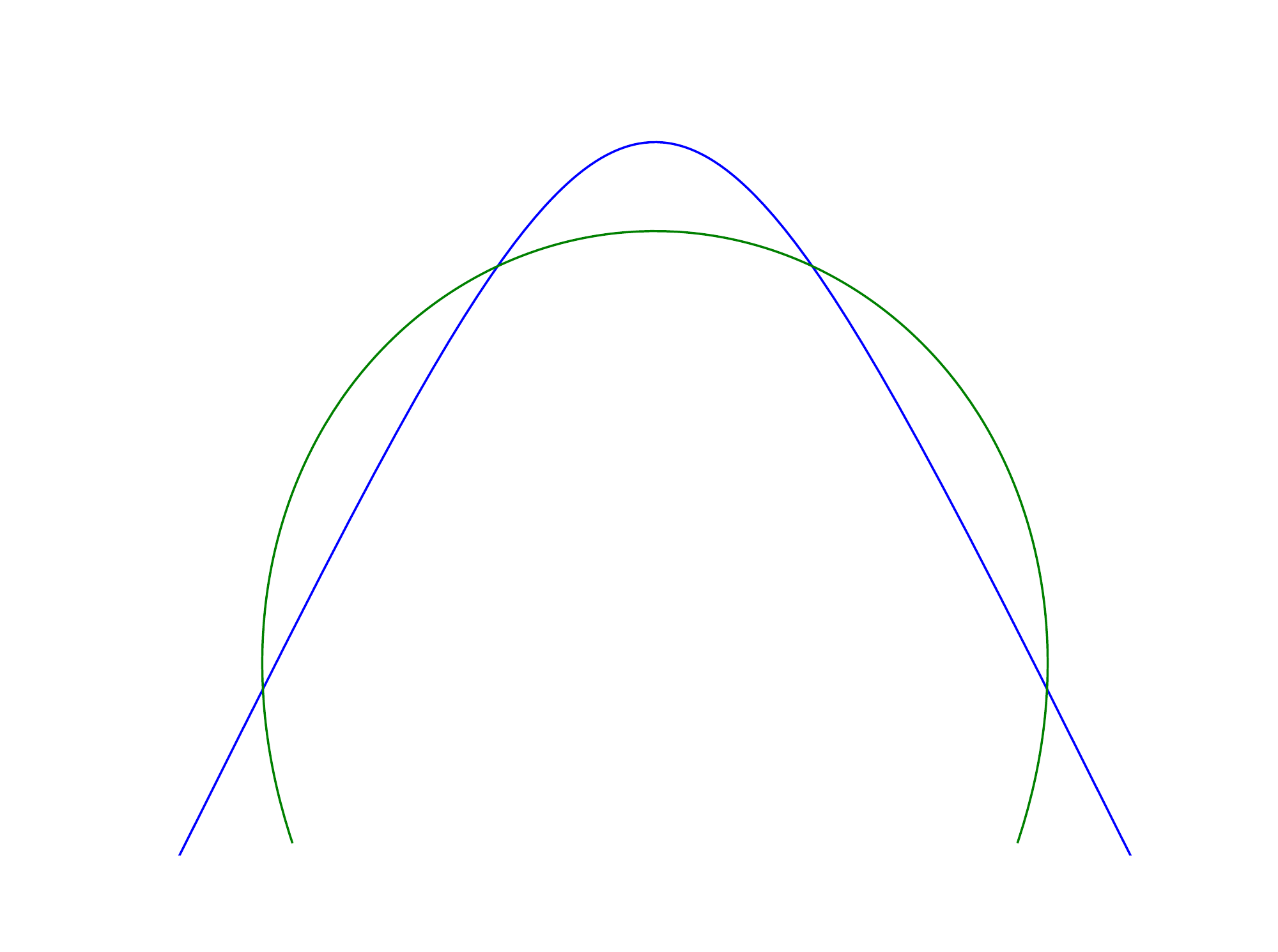}\hfil
\includegraphics[width=\unitlength,trim=38mm 30mm 38mm 30mm,clip]{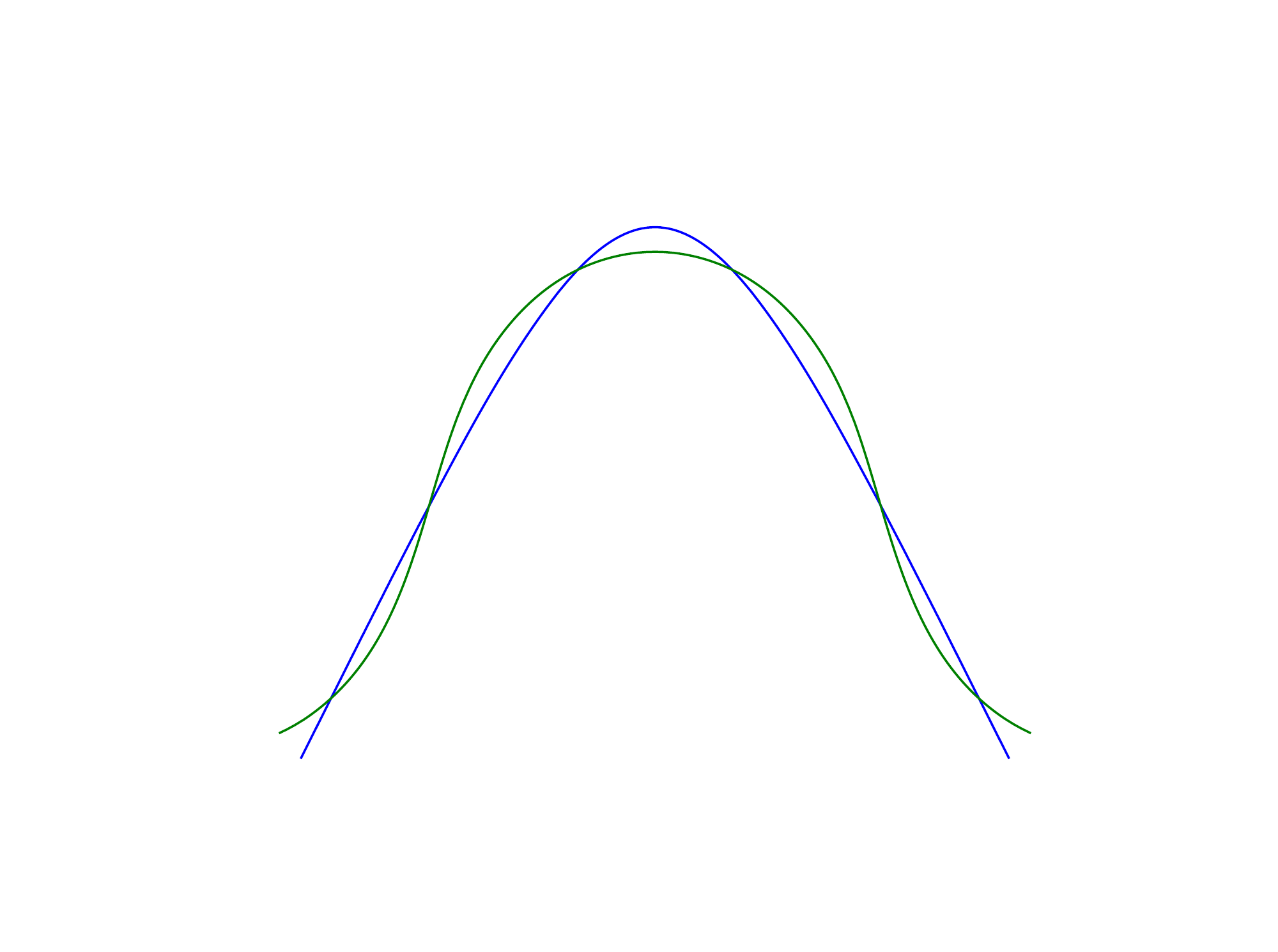}
\caption{
The blue curve is to be approximated by an elastica segment.
The green curves show the result of IPOPT optimization with different initial guesses.
In the third case the optimization terminated before an extremum was reached.}
\label{fig:init_guesses}
\end{figure}

\noindent \textbf{Outline of this article:}
We first make use of the
analytic solutions, which were derived in 1880 by L.~Saalsch\"utz \cite{saalschutz}, to
 give a representation of the space of elastic curve segments. The space
depends real analytically on seven control parameters, and one can therefore, in principle, use gradient driven optimization software to
achieve an approximation. The challenge is that the result of this non-convex nonlinear optimization depends very much on the initial guess,
as illustrated in  Figure~\ref{fig:init_guesses}.

The optimization approach
can only deliver a useable approximation algorithm if a means of choosing a good \emph{initial} elastic curve is found.
Our solution is to use the fact that the curvature function for an arc-length parameterized elastica
is \emph{affine} in a certain direction.  We first show, in Section \ref{sec:elastic-parameters},  how to use
this fact to recover the seven control parameters -- in a numerically stable manner --
 from a given elastic curve segment.   We can then (Section \ref{sec:approximating}) apply essentially
the same procedure to an \emph{arbitrary} $C^2$ plane curve
 to obtain a canonical elastic curve segment that, with respect to the global geometric characteristics selected for
in the procedure, is close to the given curve.

The canonical initial guess can then be taken as input for an optimization package to compute a good approximating
elastic curve for a given curve, provided a good solution exists.
  We illustrate the results of the overall algorithm on a sample of
B\'ezier curves -- some of which are close to elastica and some of which are not --
with endpoints free (Figure \ref{fig:opt1}) and fixed (Figure \ref{fig:opt2}).

Finally, in Section \ref{sec:conclusion}, we discuss applications of this method to the problem of approximating
curves by \emph{piecewise} elastic spline curves and  ongoing work on applications in manufacturing.


\section{Planar Elastica}
\subsection{Euler-Lagrange equation}   \label{EL-section}
Here is a brief sketch of the equations defining planar elastica.
More details, background and references can be found in \cite{Levienthesis}.
Let $\y\fr[0,\ell]\to\R^2$ be a plane curve segment parameterized by arclength. Let $\th(s)$ denote the \emph{tangent angle}, defined by the equation $\dot\y(s)=(\cos\th(s),\sin\th(s))$. Then a curve segment of length $\ell$ starting at $(x_0,y_0)$ and ending at $(x_\ell,y_\ell)$ satisfies
$x_\ell=x_0+\int_0^\ell\cos\th \ds$ and $y_\ell=y_0+\int_0^\ell\sin\th \ds$.

Let $\kappa$ denote the curvature $\dot\theta(s)$.
 An \emph{elastica} is a minimizer, among curves with the same endpoints and end tangents, of the \emph{bending energy}
$\tfrac12\int_0^\ell\k(s)^2\ds$.
Suppose $\y$ is an elastica from $(x_0,y_0)$ to $(x_\ell,y_\ell)$ with angle function $\th$, and
consider the perturbed curve $\y_t$  with angle function $\th_t(s)=\th(s)+t\psi(s)$, where $\psi$ is a differentiable function with $\psi(0)=\psi(\ell)=0$.
Applying the method of Lagrange multipliers to the bending energy, we set:
\[
\bE(\y)=\tfrac12\int_0^\ell\left(\diff{\th}s\right)^2 \ds+\l_1\left(x_0+\int_0^\ell\cos\th \ds-x_\ell\right)+\l_2\left(y_0+\int_0^\ell\sin\th \ds-y_\ell\right),
\]
and we require that
\[
0 = \diff{\bE(\y_t)}t\Big|_{t=0}
  =-\int_0^\ell\psi\left(\frac{\dd^2\th}{\ds^2} +\l_1\sin\th -\l_2\cos\th\right) \ds.
\]
Since $\psi$ was arbitrary, it follows that $\th$
satisfies the Euler-Lagrange equation
\begin{align} \label{eq:elastic-curvature-deriv-eq}
\frac{\dd^2\th}{\ds^2} +\l_1\sin\th -\l_2\cos\th=0.
\end{align}
Setting $(\lambda_1, \lambda_2) = \lambda ( \cos \phi, \sin \phi)$, with $\lambda \geq 0$, this
becomes
$\ddot \th  +\l\sin(\th-\phi)=0$.
Note that $\lambda = 0$ if and only if $\kappa$ is constant, i.e.~the curve $\y$ is either
a straight line segment or a piece of a circle.  If $\lambda \neq 0$, set
$\tilde\y(s)=\sqrt{\lambda} R_{-\phi} \y (  s/ \sqrt{\lambda} )$,
where  $R_\phi$ is the rotation by angle $\phi$.
Then $\tilde \y$ is also an elastica with tangent angle $\tilde \th(s)= \theta(s/\sqrt\lambda)-\phi$ satisfying the normalized
pendulum equation $\tilde \th ^{\prime \prime}=-\sin \tilde \th$.
In summary:
\begin{theorem}
Up to a scaling and rotation of the ambient space, all arclength parameterized elastica
 $\y\fr[0,1] \to \R^2$, with non-constant curvature $\kappa$, can be expressed as:
$\y(s)= \y(0) + \int_0^s \left( \cos \th(t), \,\, \sin \th(t)\right)\dt$,
where
\begin{align}\label{eq:elas-diff}
\ddot\th=-\sin\th.
\end{align}
\end{theorem}


\subsection{Parameterizations of the space of elastica}
We now want to define some suitable control parameters to describe an arbitrary
elastic curve segment. Essentially, the parameters need to specify which
solution to \eqref{eq:elas-diff} is involved, the start and endpoints on the solution curve in question, and
a rotation and scaling.

The elastic curves can be expressed in closed form via the elliptic functions (see Appendix~\ref{sec:app_elliptic}). The formulas can be found in Love \cite{Love}.
There are two classes of elastica: curves with inflection points (i.e. points where $\dot\th=0$) and curves without inflections.

\subsubsection*{Basic elastica}
\begin{figure}[h!] \centering
\includegraphics[width=\textwidth,trim=10mm 30mm 30mm 30mm,clip]{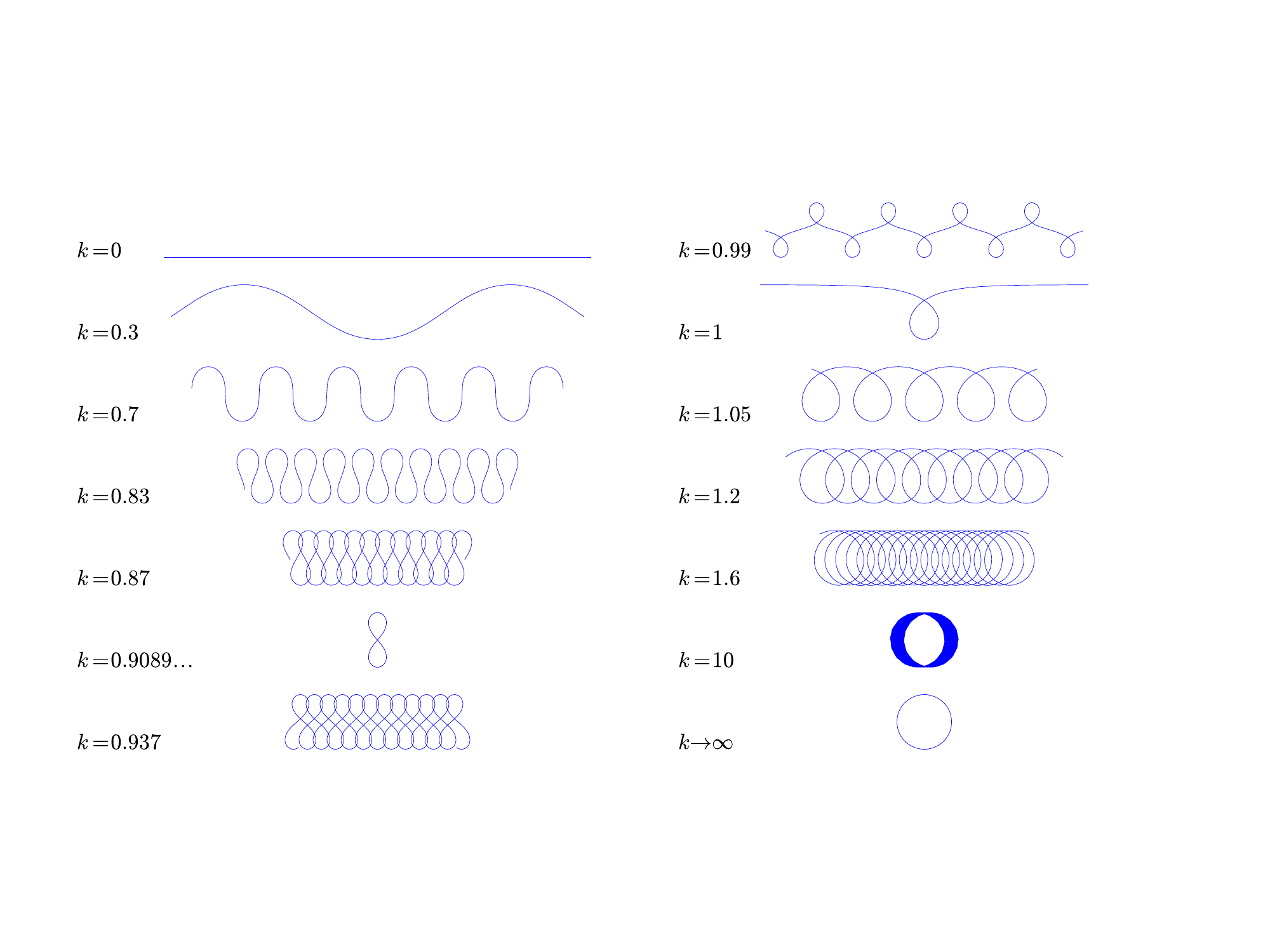}
\caption{Elastica for different values of $k$. The curves are scaled to a uniform ``height".}\label{fig:elas_pic}
\end{figure}

The solution to \eqref{eq:elas-diff} starting at $(0,0)$ with initial angle $\th(0)=0$ and $\dot\th(0)\geq0$ is
\[
\infl_k(s)=\left( 2E(s,k)-s , \,\,  2k(1-\cn(s,k)) \right),
\]
where $k=\dot\th(0)/2$. For $k\in[0,1 )$ we get inflectional elastica;  for $k\geq1$,
we use the \emph{extended} elliptic functions defined in Appendix~\ref{sec:app_elliptic} to obtain
 elastic curves without inflections. We reserve the name $\infl_k$ for these basic elastica. Figure \ref{fig:elas_pic} shows elastic curves for different values of $k$. All elastica are obtained by scaling and rotating these curves.
The periodicity of the curves is given by:
\[ \infl_k(s+4K)=\infl_k(s)+\left( 2E(4K)-4K , \, \, 0 \right),
\]
where $K$ the quarter-period defined in Appendix~\ref{sec:app_elliptic}.

\subsubsection*{General elastica}
To parameterize an arbitrary segment of an elastica, we can choose a segment of a basic elastica by choosing $k$, a starting point $\stp$ and an endpoint
$\stp + \ell$, where $\ell \in \R \setminus\{0\}$. It will be convenient to remove the dependence on $\stp$ and $\ell$ from the domain and include them in the parameterization. We parameterize the elastica segments on the unit interval, setting $s=\stp+\ell t$, with $t\in[0,1]$. The new curve parameter $t$ is not unit speed.
Finally, any elastica segment can be obtained by introducing a scaling factor $\w>0$, a rotation by an angle $\phi\in\left(-\pi,\pi\right]$ and translation by a vector $(x_0,y_0)$. We thus have a standard elastic segment parameterization
\[
\y_{(k,\stp,\ell,\w,\phi,x_0,y_0)}(t)=\w R_\phi\infl_k(\stp+\ell t)+
      (x_0 , \,  y_0 ), \quad t\in[0,1].
\]
It depends on seven control parameters,  but we will usually omit the subscript. Such a curve
 has constant speed $|\ell|\w$ and length  $L=|\ell|\w$.

\begin{note}If $\ell<0$, the orientation of the curve is changed. In the inflectional case, the elastica with opposite orientation can be obtained by a rotation by $\pi$. In this case we may therefore assume $\ell>0$ without loss of generality. For elastica \emph{without} inflections, however, a segment with $\ell<0$ cannot be described as a segment with $\ell>0$.  One can instead reverse the direction of the parameterization for that case, and
so all cases can be handled with the assumption $\ell>0$. \end{note}

For any elastic curve  $\y$ of the above type, the curve $\tilde\y(t)=\y(\frac t{\ell\w})$ is unit speed.
Letting $\th$ and $\tilde\th$ denote the angle functions of $\infl_k$ and $\tilde\y$, respectively, we have $\tilde\th(t)=\th(\stp+\frac t\w)+\phi$, and thus
 \[
\tilde\th''(t)=\tfrac1{\w^2}\ddot\th(\stp+\tfrac t\w)=-\tfrac1{\w^2}\sin\th\left(\stp+\tfrac t\w\right)=-\tfrac1{\w^2}\sin\left(\tilde\th(t)-\phi\right),
\]
so the angle function for $\tilde\y$ satisfies \eqref{eq:elastic-curvature-deriv-eq} with
 \begin{align}\label{eq:lambdas} ( \lambda_1 , \,  \lambda_2) =\tfrac1{\w^2}(\cos\phi, \, \sin\phi ).
\end{align}

 We will also use the fact that the curvature for the elastica $\y_{(k,\stp,\ell,\w,\phi,x_0,y_0)}$ is
\begin{align}\label{eq:kappa-infl} \kappa(t)=\tfrac{2k}\w\cn(\stp+\ell t).\end{align}


\section{Finding the control parameters of an elastic curve segment}\label{sec:elastic-parameters}
We first describe a way to calculate numerically the control parameters of a given planar elastic curve segment.
In the next section the same recipe will be applied to an arbitrary planar curve to obtain a canonical
 first guess for an approximating elastic curve. The main idea is to exploit the fact that the curvature of an
 elastica is an affine function of the distance along a special direction.

Let $\bbx:[a,b]\to\R^2$ be an elastic curve parameterized by arclength. As for any planar curve,
we can write the tangent and the normal as
$\bbt=(\cos\theta,\sin\theta)$, $\bbn=(-\sin\theta,\cos\theta)$, and
we have the Frenet-Serret equations
\begin{align*}
  \frac{\dd\bbt}{\ds}=\frac{\dtheta}{\ds}\,\bbn=\kappa\,\bbn\,, \quad
  \frac{\dd^2\bbt}{\ds^2}=\frac{\dd^2\theta}{\ds^2}\,\bbn-\kappa^2\,\bbt\,.
\end{align*}
The tangent angle $\theta$ must  satisfy the Euler-Lagrange equation \eqref{eq:elastic-curvature-deriv-eq}
for some Lagrangian multipliers $\lambda_1,\lambda_2$ to be found.

Let $u$ denote the projection of $\bbx$ onto the line spanned by
$(\lambda_2,-\lambda_1)$, i.e.,
\begin{align*}
  u&=\frac{1}{\lambda} ( \lambda_2 , \,  -\lambda_1 ) \, \cdot \, ( x, \,  y )
  =\frac{\lambda_2\,x-\lambda_1\,y}{\lambda},
\end{align*} where $\l=\|(\lambda_1,\lambda_2)\|=S^{-2}$.
Setting $\phi=0$ in \eqref{eq:lambdas}, we find that the vector $(\lambda_2, -\lambda_1)$ points in the
 downward direction in Figure \ref{fig:elas_pic}.  It follows that $u$ is bounded and periodic in $s$.
Moreover, we can write the Euler-Lagrange equation as $\ddot \theta =\lambda\, \dot u$, so we have
\begin{equation}
  \label{eq:alpha-def}
  \kappa=\frac{\dtheta}{\ds}
  =\lambda\,u+\alpha=\lambda_2\,x-\lambda_1\,y+\alpha\,,
\end{equation}
which is to say that the curvature is an affine function of $u$.

In order to find $\lambda_1$, $\lambda_2$ and $\alpha$ in a numerically stable manner, we solve the above equation in the least squares sense, i.e., we
consider the quadratic minimization problem
\begin{equation*}
  \minimize_{\lambda_1,\lambda_2,\alpha}\int_{a}^{b}\left(
    \kappa+\lambda_1\,y-\lambda_2\,x-\alpha
  \right)^2\,\ds\, ,
\end{equation*}
which leads to the following linear system
\begin{equation}
  \label{eq:lsq-curvature-eq}
  \begin{pmatrix}
    \int_{a}^{b}y^2\,\ds
    &-\int_{a}^{b}x\,y\,\ds
    &-\int_{a}^{b}y\,\ds
    \\
    -\int_{a}^{b}x\,y\,\ds
    &\int_{a}^{b}x^2\,\ds
    &\int_{a}^{b}x\,\ds
    \\
    -\int_{a}^{b}y\,\ds
    &\int_{a}^{b}x\,\ds
    &\int_{a}^{b}1\,\ds
  \end{pmatrix}
  \begin{pmatrix}
    \lambda_1\\\lambda_2\\\alpha
  \end{pmatrix}
  =
  \begin{pmatrix}
    -\int_{a}^{b}y\,\kappa\,\ds
    \\
    \int_{a}^{b}x\kappa\,\ds
    \\
    \int_{a}^{b}\kappa\,\ds
  \end{pmatrix}
  \,.
\end{equation}

\begin{figure}[h]
  \centering
  \unitlength=.7\textwidth
  \includegraphics[width=\unitlength,trim=30mm 20mm 30mm 20mm,clip]{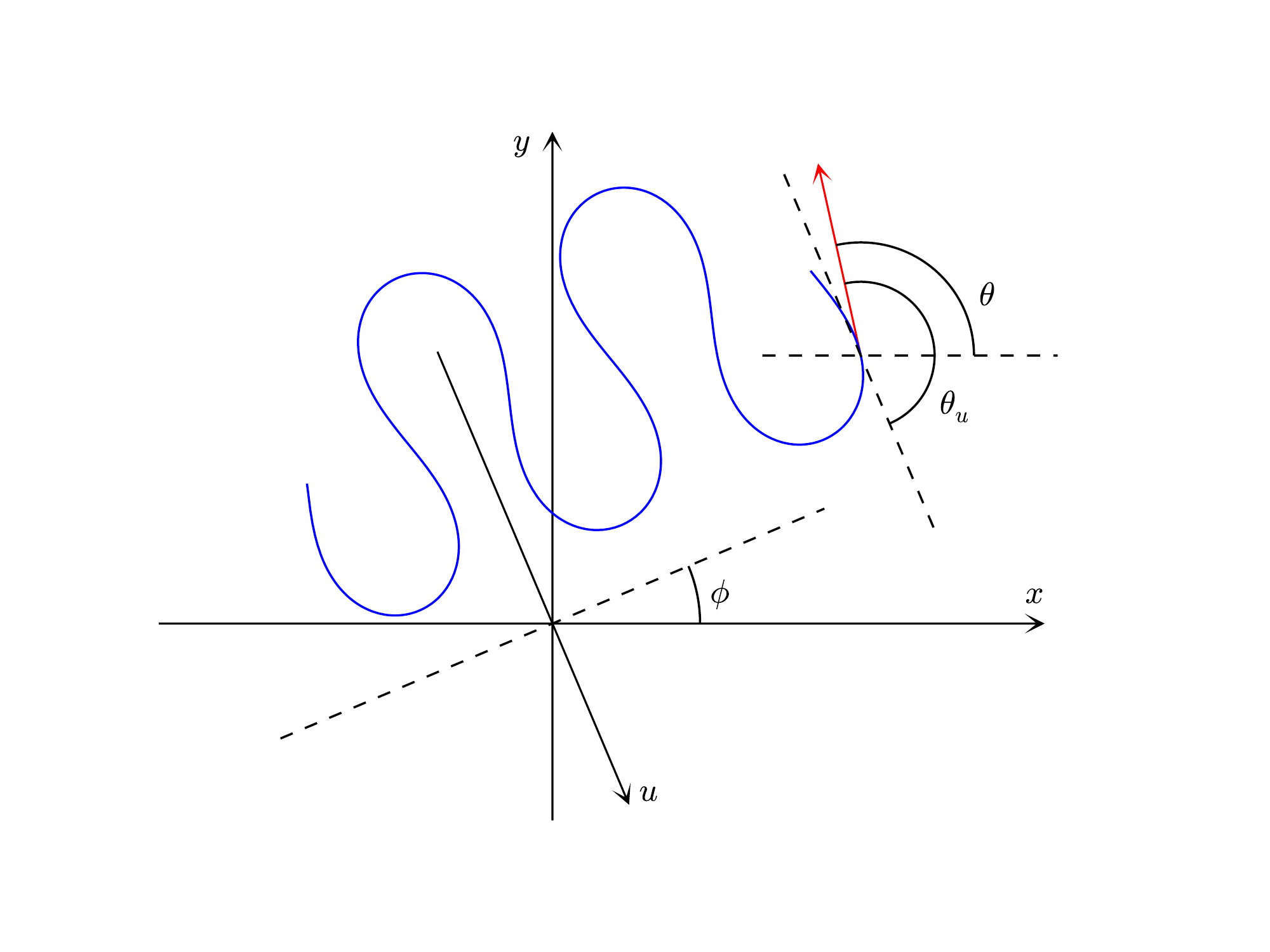}
  \caption{A segment of an elastica in the plane, showing $\phi$,  $\theta$, the $u$-axis and the angle $\theta_u$.}
  \label{fig:phi-theta}
\end{figure}

Let $\theta_u$ denote the angle between the tangent vector
$\bbt$ and the $u$-axis (see Figure~\ref{fig:phi-theta}). We have
\begin{align}
  \label{eq:theta_e}\begin{split}
    \cos\theta_u&=\frac{1}{\lambda}( \lambda_2, \, \, -\lambda_1 ) \cdot \bbt
    =\frac{1}{\lambda}( \lambda_2, \, \, -\lambda_1) \cdot\frac{\dd\bbx}{\ds}
    =\frac{\du}{\ds}\,,\\
    \sin\theta_u&=\frac{1}{\lambda} ( \lambda_1, \, \, \lambda_2) \cdot\bbt\,,\end{split}
\end{align}
and hence
\begin{equation*}
  \frac{\dd\sin\theta_u}{\du}
  =\frac{\ds}{\du}\,\frac{\dd\sin\theta_u}{\ds}
  =\frac{1}{\cos\theta_u}\,\cos\theta_u\,\frac{\dtheta_u}{\ds}
  =\kappa=\lambda\,u+\alpha\,,
\end{equation*}
or equivalently
\begin{equation}
  \label{eq:beta-def}
  P(u) := \sin\theta_u=\frac12\lambda\,u^2+\alpha\,u+\beta \,.
\end{equation}
We solve this equation with respect to $\beta$ in the least squares
sense and obtain
\begin{equation}
  \label{eq:beta-eq}
  \beta=\frac{1}{L}
  \int_{a}^{b}\left(\sin\theta_u-\frac12\lambda\,u^2-\alpha\,u\right)\,\ds\, ,
\end{equation} where $L=b-a$ is the length of the curve $\bbx$.

For an elastica $\bbx(s)=SR_\phi \zeta_k(s/S)+(x_0,y_0)$ we have
$(\lambda_1, \lambda_2) = S^{-2}(\cos \phi, \sin \phi)$ and $\kappa(s)=(2k/S)\cn(s/S)$.
 Substituting these into the definitions
$u=(\lambda_2 x - \lambda_1 y)/\lambda$, $\alpha = \kappa - \lambda u$ and
$\sin \theta_u =  (1/\lambda)(\lambda_1, \lambda_2) \cdot {\bf t}$ we have
\begin{eqnarray}
u &=& -2Sk(1-\cn(s/S)) + x_0 \sin \phi - y_0 \cos \phi, \nonumber \\
\sin \th_u &= & 2\dn^2(s/S)-1,  \nonumber\\
  \alpha & = & 2 k S^{-1} \cn(s/S) - \lambda u =
	 2k/S+(y_0\cos\phi-x_0\sin\phi)/S^{2}   \label{eq:alpha-parameters}
\end{eqnarray}

Then the equation $\beta = \sin \theta_u - \lambda u^2/2 - \alpha u$ becomes:
 \begin{align}\label{eq:beta-parameters}\beta= 1+\frac{x_0\sin\phi-y_0\cos\phi}{2\w^2}(x_0\sin\phi-y_0\cos\phi-4k\w).\end{align}

It follows from \eqref{eq:beta-def} that all points on the elastica correspond to $u$-values where the
value of the polynomial $P$ is between $-1$ and $1$, and hence \[ u\in\left[\frac{-\a-\delta_-}{\l},\frac{-\a+\delta_-}{\l}\right],\] where $\delta_-=\sqrt{\a^2-2\l(\b-1)}$.

\begin{figure}[here]
\centering
\includegraphics[totalheight=50mm,trim=30mm 20mm 40mm 20mm,clip]{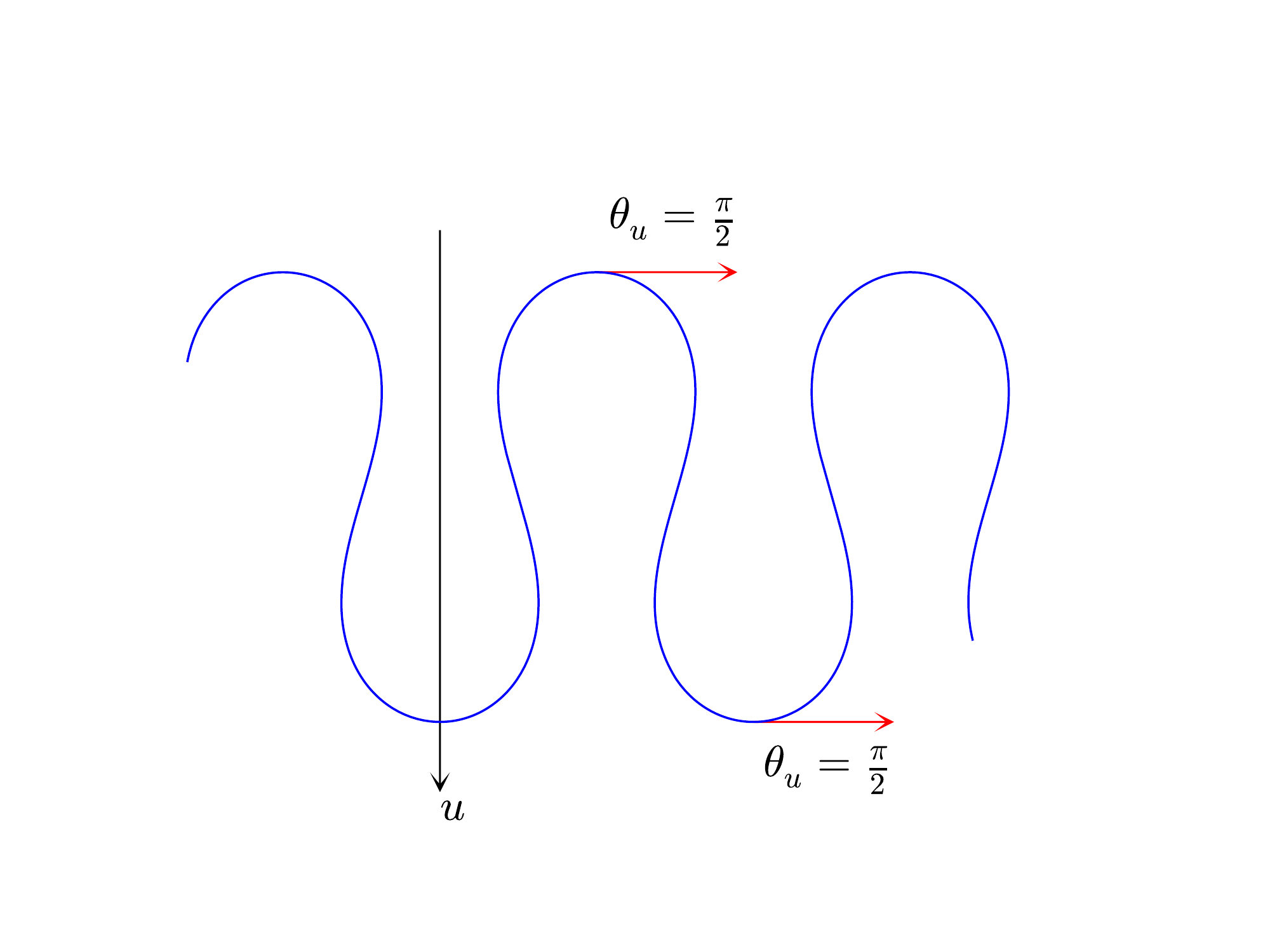} \hfil
\includegraphics[totalheight=50mm,trim=30mm 16mm 20mm 33mm,clip]{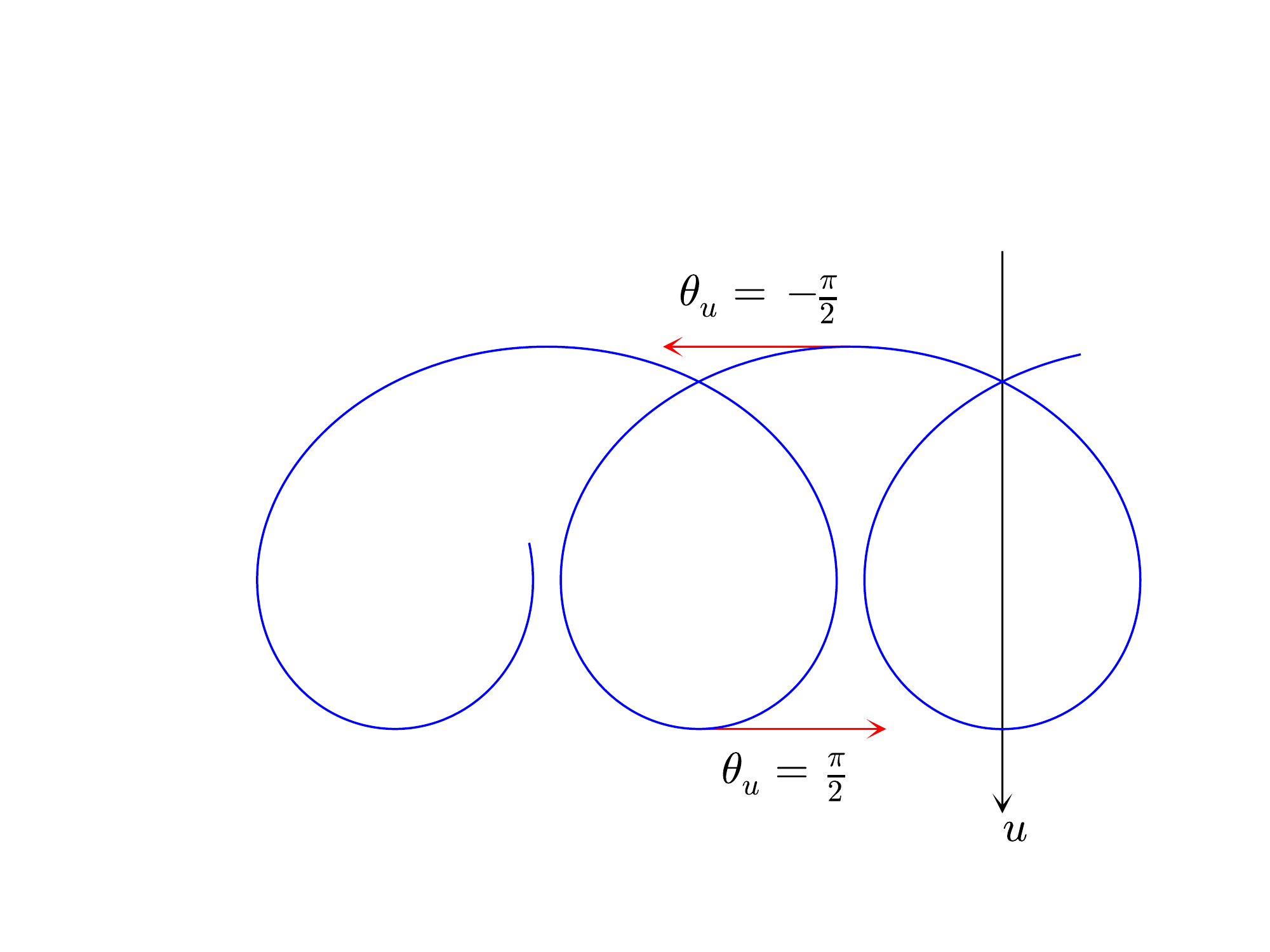}
\caption{The inflectional elastica have points where $\sin\th_u=1$, but not $-1$. For the non-inflectional elastica $\sin\th_u$ takes both the values $\pm1$.}
\label{fig:theta-u}
\end{figure}

If the elastica has an inflection, there must be some $u_*$, such that $\k=\l u_*+\a=0$, but this means that $u_*$ is a minimizer for $P(u)$, and thus its minimum must lie in $[-1,1]$. Moreover, the inflectional elastica has points where $\sin\th_u=1$ (which happens twice per period), but no points where $\sin\th_u=-1$ (see Figure~\ref{fig:theta-u}). Hence $u$ runs through all of the interval where $P(u)$ is less that $1$; in other words, $u_{\min}$ and $u_{\max}$ are exactly the endpoints of the above interval.

\begin{figure}[h]
  \centering
  \unitlength=.46\textwidth
  \begin{picture}(1,0.75)
  \includegraphics[width=\unitlength]{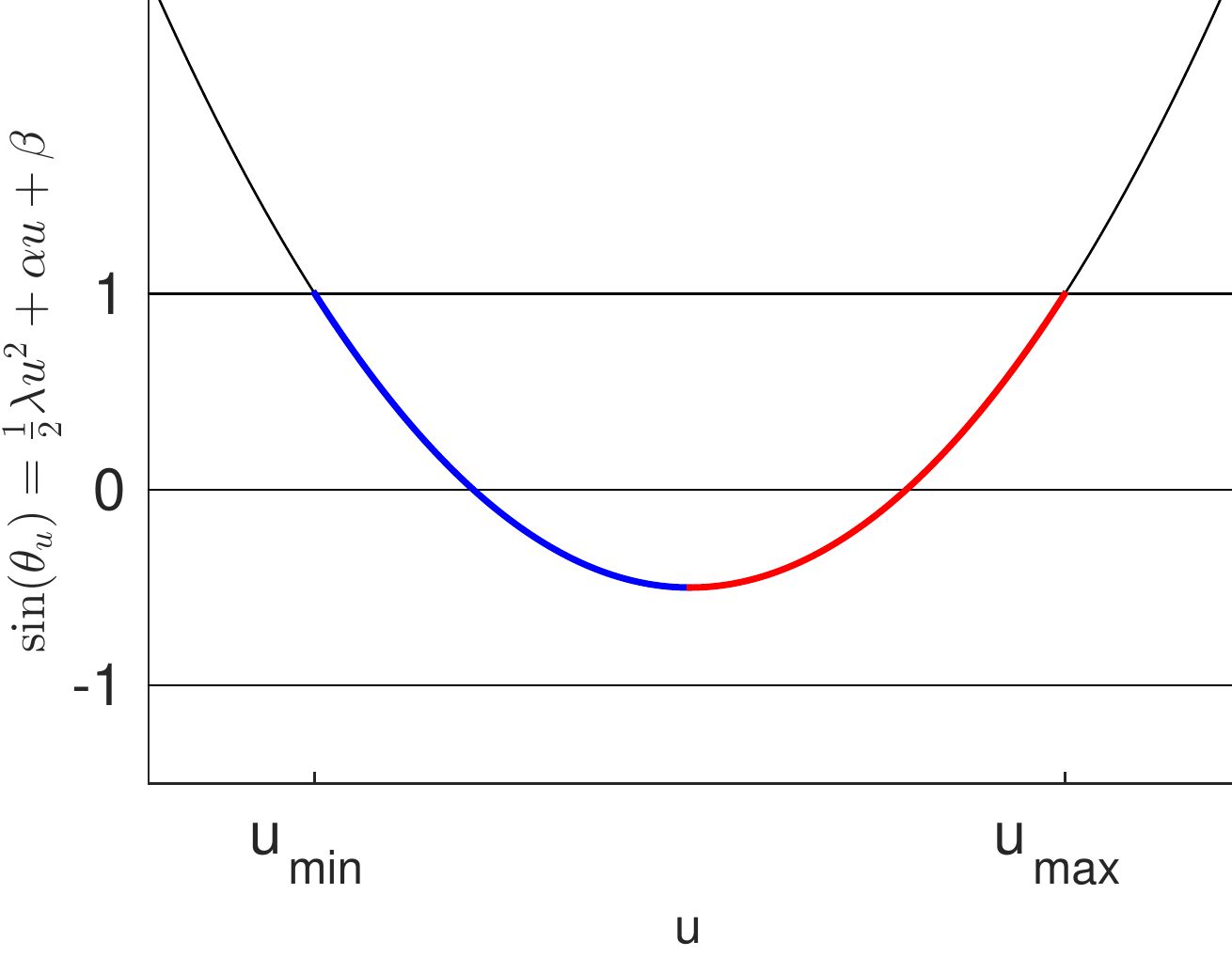}
  \end{picture}\quad \quad
  \begin{picture}(1,0.75)
  \small
  \put(0,0){\includegraphics[width=\unitlength]{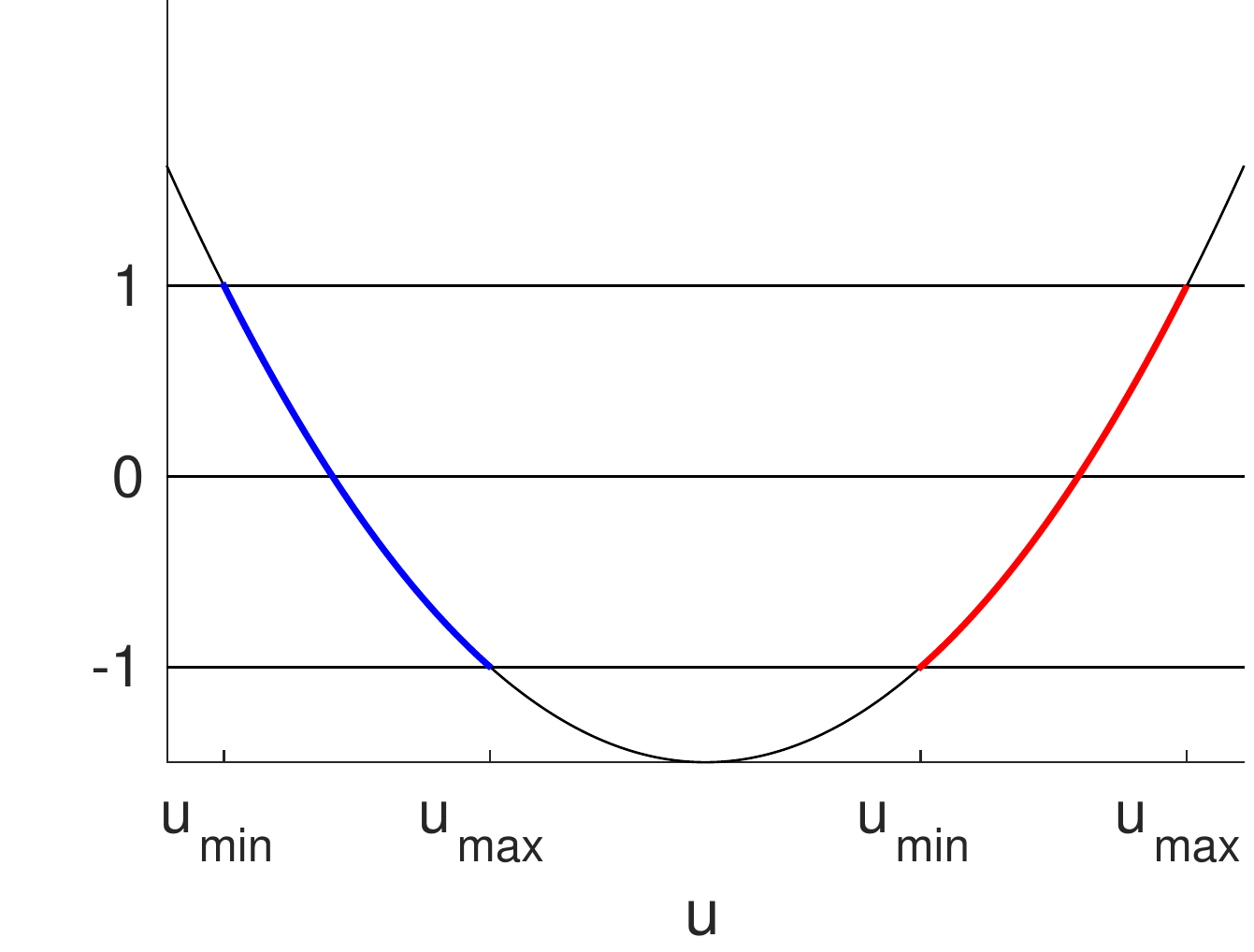}}
  \put(0.25,0.7){$\ell<0$}
  \put(0.75,0.70){$\ell>0$}
  \end{picture}\hfil
  \caption{The parabola \eqref{eq:beta-def}. To the left in the case
    of an elastica with inflection points, to the right without. The
    blue and red part corresponds to points with negative and positive
    curvature, respectively.}
  \label{fig:parabola}
\end{figure}

For the elastica without inflections, the tangent makes full rotations, so $\sin\th_u$ takes both of the values $\pm1$. Hence, in this case,
 we must have
\begin{align*} [u_{\min},u_{\max}]=\left[\frac{-\a-\delta_-}{\l},\frac{-\a-\delta_+}{\l}\right] \quad \textrm{ or }\quad [u_{\min},u_{\max}]=\left[\frac{-\a+\delta_+}{\l},\frac{-\a+\delta_-}{\l}\right],\end{align*}
where $\delta_+=\sqrt{\a^2-2\l(\b+1)}$;  these are the two cases corresponding
 to $\ell<0$ and $\ell>0$, respectively.  We can thus determine whether the elastica has inflection points based on whether the minimum for the polynomial $P(u)$ is smaller or greater than $-1$, see Figure~\ref{fig:parabola}. In fact, from \eqref{eq:alpha-parameters} and \eqref{eq:beta-parameters} we have \begin{align*}
\a^2-2\l(\b-1)=\frac{4k^2}{\w^2},
\end{align*} or equivalently
\begin{equation}
  \label{eq:k-type1}
  k=\frac{\sqrt{\a^2-2\l(\b-1)}}{2\sqrt{\l}}\,,
\end{equation}
so we can find $\w$, $\phi$ and $k$ from $\l_1$, $\l_2$, $\a$ and $\beta$.
\begin{note} The above formula also holds if $\ell<0$ and so does the expression for $\b$ in control parameters. The expressions for $\kappa$, $\l_1$, $\l_2$ and $\a$ simply change sign in this case.\end{note}

We still need to recover $\stp$ and $\ell$. We have \[u=-2k\w\left(1-\cn\left(\stp+\tfrac t\w\right)\right)+x_0\sin\phi -y_0\cos\phi,\]
and since
\[ u_{\max}=\frac{-\a+\delta_-}{\l}=x_0\sin(\phi)-y_0\cos(\phi), \]
we get
\begin{equation*}
  \Delta(u)=u_{\max}-u=2k\w\left(1-\cn s\right)\,,
\end{equation*}
so
\begin{equation}
\label{eq:cn_umax}
\cn(s,k)=1-\frac{\Delta(u)}{2k\w}\,.
\end{equation}
If we consider the unbounded complete elastica, then $u$ oscillates
between $u_{\min}$ and $u_{\max}$ and we can divide the elastica into
segments where $u$ is monotone, each with length equal to a half
period $2K\w$.

We first consider the case of an elastica with inflection points (i.e. $k<1$). Here we have $\cn(s,k)=\cos\left(\am(s,k)\right)$.
If the start point
$\bbx_0=\bbx(a)$ is on segment number 1 and $u$ is decreasing here,
then
\begin{equation*}
  \am(s_0,k)=\arccos\left(
    1-\frac{\Delta(u_0)}{2k\w}
  \right)\,,
\end{equation*}
and if the end point $\bbx_1=\bbx(b)$ is on segment number $n$, then
\begin{equation*}
  \am(s_1,k)=
  \begin{cases}
    (n-1)\,\pi+\arccos\left(
    1-\frac{\Delta(u_1)}{2k\w}
  \right)\,,
  &\text{if $n$ is odd,}
  \\
  n\,\pi-\arccos\left(
    1-\frac{\Delta(u_1)}{2k\w}
  \right)\,,
  &\text{if $n$ is even.}
  \end{cases}
\end{equation*}
If $u$ is increasing on segment number 1, then
\begin{equation*}
  \am(s_0,k)=2\pi-\arccos\left(
    1-\frac{\Delta(u_0)}{2k\w}
  \right)\,,
\end{equation*}
and
\begin{equation*}
  \am(s_1,k)=
  \begin{cases}
    (n+1)\,\pi-\arccos\left(
    1-\frac{\Delta(u_1)}{2k\w}
  \right)\,,
  &\text{if $n$ is odd,}
  \\
  n\,\pi+\arccos\left(
    1-\frac{\Delta(u_1)}{2k\w}
  \right)\,,
  &\text{if $n$ is even.}
  \end{cases}
\end{equation*}
In all cases we have
\begin{equation*}
  s_i=F(\am(s_i,k),k)\,,\qquad i=0,1\,,
\end{equation*}
and $\ell=s_1-s_0$.

In the case of an elastica without inflections points (i.e. $k\geq1$) we need a little work to find $\am$. We have
\[
\sn(s,k)=\frac1k\sn\left(ks,\tfrac1k\right)=\frac1k\sin\left(\am\left(ks,\tfrac1k\right)\right)\,
\]
and 
\[
\sn(s,k)=\sqrt{1-\cn^2(s,k)}=\sqrt{\frac{\Delta(u)}{k\w}\left(1-\frac{\Delta(u)}{4k\w}\right)}\,.
\]
If $u$ is decreasing on segment 1 then
\begin{equation*}
  \label{eq:ams0-decrease-type2}
  \am\left(ks_0,\tfrac1k\right)=\arcsin\sqrt{\frac{\Delta(u_0)}{\w}\left(k-\frac{\Delta(u_0)}{4\w}\right)}\,,
\end{equation*}
and if we have $n$ segments
\begin{equation*}
  \label{eq:ams1-decrease-type2}
  \am\left(\frac{s_1}{k},k\right)=
  \begin{cases}
    \frac{n-1}{2}\pi+\arcsin\sqrt{\frac{\Delta(u_1)}{\w}\left(k-\frac{\Delta(u_1)}{4\w}\right)}\,,&\text{if $n$ is odd,}
    \\
    \frac{n}{2}\pi-\arcsin\sqrt{\frac{\Delta(u_1)}{\w}\left(k-\frac{\Delta(u_1)}{4\w}\right)}\,,&\text{if $n$ is even.}
  \end{cases}
\end{equation*}
If $u$ is increasing on segment 1 then
\begin{equation*}
  \label{eq:ams0-increase-type2}
  \am\left(\frac{s_0}{k},k\right)=\pi-\arcsin\sqrt{\frac{\Delta(u_0)}{\w}\left(k-\frac{\Delta(u_0)}{4\w}\right)}\,,
\end{equation*}
and if we have $n$ segments
\begin{equation*}
  \label{eq:ams1-increase-type2}
  \am\left(\frac{s_1}{k},k\right)=
  \begin{cases}
    \frac{n+1}{2}\pi-\arcsin\sqrt{\frac{\Delta(u_1)}{\w}\left(k-\frac{\Delta(u_1)}{4\w}\right)}\,,&\text{if $n$ is odd,}
    \\
    \frac{n}{2}\pi+\arcsin\sqrt{\frac{\Delta(u_1)}{\w}\left(k-\frac{\Delta(u_1)}{4\w}\right)}\,,&\text{if $n$ is even.}
  \end{cases}
\end{equation*}
Finally, we find the $s$-values using the incomplete elliptic integral
\begin{equation*}
  s_i=\tfrac1k F\left(\am\left(ks_i,\tfrac1k\right),\tfrac1k\right)\,,\qquad i=0,1\,,
\end{equation*}
and $\ell=s_1-s_0$.

\begin{note}If we have a negatively curved noninflectional elastica (i.e. $\ell<0$), we can reverse
the parameterization, find the elastica, and interchange $(s_0,s_1)$.\end{note}

We now have a scaled and rotated elastica segment, $\y_0=\y_{(k,\stp,\ell,\w,\phi,0,0)}$, and all that is left is to find the final translation $(x_0,y_0)$. This is done by solving the
equation
\begin{equation*}
  \bbx(s)=\y_0(s)+ (x_0 , \,  y_0 ),
\end{equation*}
in the least squares sense. The solution is
\begin{equation}
  \label{eq:translate}
  (x_0 , \, y_0 ) =\frac{1}{L}\int_a^b(\bbx(s)-\y_0(s))\,\ds.
\end{equation}


\section{Approximating a plane curve by a planar elastica}  \label{sec:approximating}
We are now given a curve $\bbx\fr[0,1]\to\R^2$, not necessarily elastic and
not necessarily parameterized by arclength. The arclength is given by \[s(t)=\int_0^t\|\bbx'(\tau)\|\dd\tau\,,\] and the length of the curve is $L=s(1)$. We want to approximate this curve by a piece of an elastica. We do this by minimizing a suitable
distance, such as the $L^2$, $H^1$, or $H^2$ distance, over the control parameters $\mathbf{p}=(k,\stp,\ell,\w,\phi,x_0,y_0)$. In the case of $L^2$ the problem is
\[\minimize_{k,\stp,\ell,\w,\phi,x_0,y_0}\mathcal{F}(k,\stp,\ell,\w,\phi,x_0,y_0),\] where
\begin{align*}
\mathcal{F}(\mathbf{p})=\frac{1}{2}\int_0^1\left\|\y_{\mathbf{p}}\left(\frac{s(t)}{L}\right)-\bbx(t)\right\|^2\|\bbx'(t)\| \,\dt\,.
\end{align*}

If we want the elastic curve to satisfy further conditions, such as having the same endpoints and/or end tangents as the original curve, we can include these in the optimization problem.

For the optimization, we have used the gradient driven tool IPOPT \cite{ipopt}, so we need the first and second order partial derivatives of $\mathcal{F}$ with respect to the control parameters, which are
\begin{align*} \frac{\partial\mathcal{F}}{\partial p_i}
 = &\int_0^1\left(\y_{\mathbf{p}}\left(\frac{s(t)}{L}\right)-\bbx(t)\right)\cdot\tfrac{\partial\y_{\mathbf{p}}}{\partial p_i}\left(\frac{s(t)}{L}\right)\|\bbx'(t)\|\dt\,, \\ \frac{\partial^2\mathcal{F}}{\partial p_i\partial p_j}
= &\int_0^1\left(\tfrac{\partial\y_{\mathbf{p}}}{\partial p_i}\left(\frac{s(t)}{L}\right)
 \cdot\tfrac{\partial\y_{\mathbf{p}}}{\partial p_j}\left(\frac{s(t)}{L}\right)   \right. \\
  &\left. +  \left(\y_{\mathbf{p}}\left(\frac{s(t)}{L}\right) -\bbx(t)\right)\cdot\tfrac{\partial^2\y_{\mathbf{p}}}{\partial p_i\partial p_j}
\left(\frac{s(t)}{L}\right) \right)\|\bbx'(t)\|\dt\,.\end{align*}
See Appendix~\ref{sec:app_derivatives} for a list of specific derivatives.

The optimization problem is non convex and the result depends on the initial guess (see Figure~\ref{fig:init_guesses}).
A canonical geometrically plausible guess is obtained from a generalization of
 the procedure of Section~\ref{sec:elastic-parameters} to the case of an arbitrary input curve, which we will now
describe.

We find $\lambda_1,\lambda_2,\alpha$ as before, by solving
\eqref{eq:lsq-curvature-eq}, noting that
  $\int_{s(0)}^{s(1)}f\,\ds=\int_0^1f(t)\,\frac{\ds}{\dt}\,\dt$.
This gives us the scaling and rotation of the elastica. We can judge
the success by calculating the normalized residual
\begin{equation*}
  R_1=\sqrt{\int_0^1\left(\kappa(t)+\lambda_1\,y(t)-\lambda_2\,x(t)-\alpha\right)^2
    \,\frac{\ds}{\dt}\,\dt}
  \Bigg/\sqrt{\int_0^1\kappa^2(t)\,\frac{\ds}{\dt}\,\dt}\,.
\end{equation*}
Similarly $\beta$ can be found by \eqref{eq:beta-eq}, where
$\sin\theta_u$ is given by \eqref{eq:theta_e}, and we
can calculate the normalized residual
\begin{equation*}
  R_2=\sqrt{\frac{1}{L}\int_0^1\left(
      \sin\theta_u(t)-\tfrac12\lambda\,u^2(t)-\alpha\,u(t)-\beta
    \right)^2\,\frac{\ds}{\dt}\,\dt}\,.
\end{equation*}

\begin{note}Another possibility is to forget that we know $\lambda$ and $\alpha$ and solve
\eqref{eq:beta-def} with respect to $\lambda$, $\alpha$, and $\beta$, but in the few cases we tried this, the results got worse.
\end{note}

We know that $\sin\th_u$ takes values in $[-1,1]$, and since $\b$ is chosen to minimize the distance between $\sin\th_u$ and the polynomial $P(u)=\frac12\lambda\,u^2+\alpha\,u+\beta$, the latter must be less than $1$ for some $u$-values, so the number $\delta_-=\sqrt{\a^2-2\l(\b-1)}$ is well-defined. We can thus determine whether the elastica has inflection points
and we can determine the parameter $k$ from \eqref{eq:k-type1}.

\begin{figure}[h]
  \centering
  \includegraphics[totalheight=38mm,trim=10mm 32mm 20mm 24mm,clip]{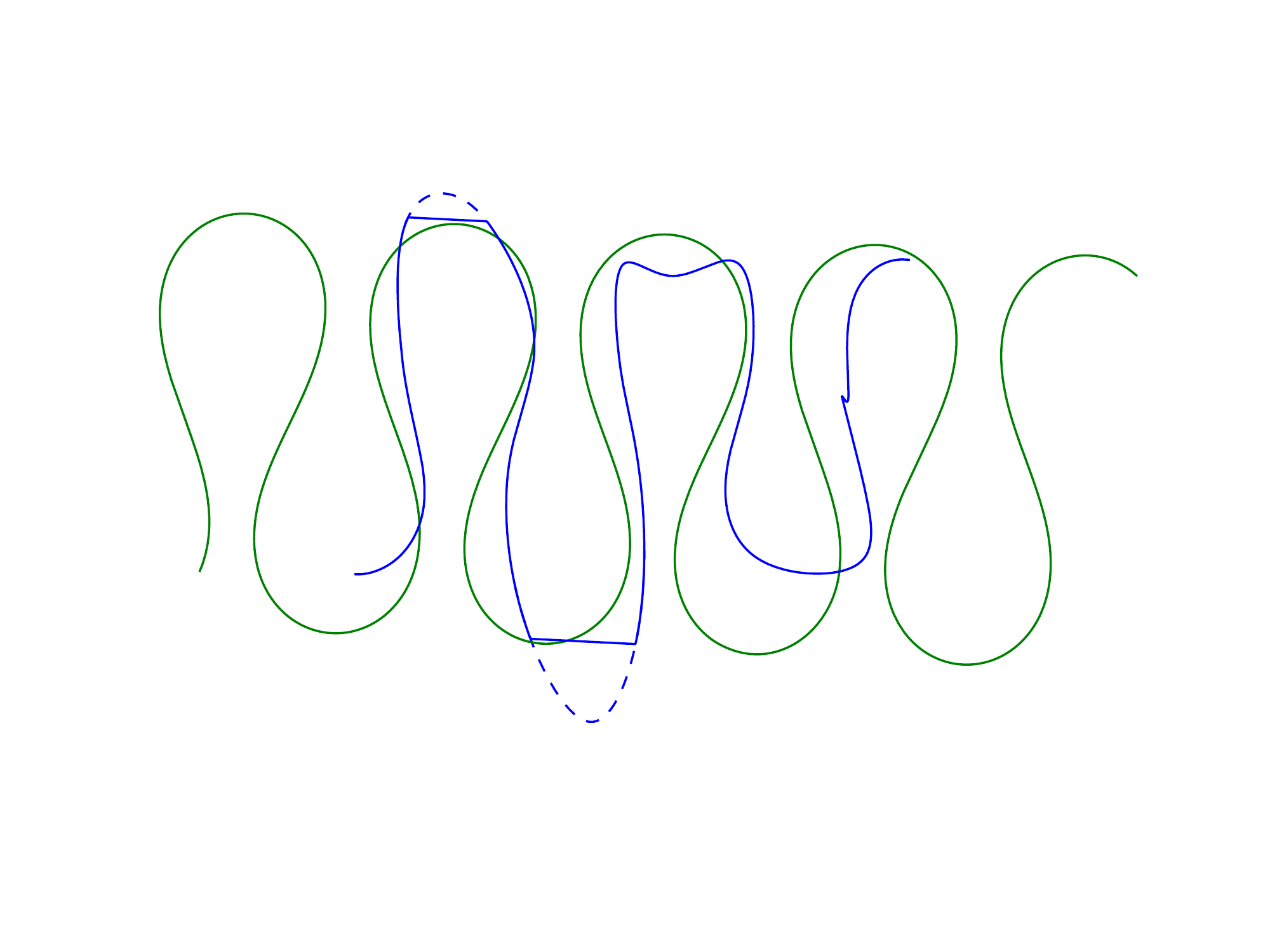}   \hfil
	\includegraphics[totalheight=38mm,trim=30mm 20mm 30mm 10mm,clip]{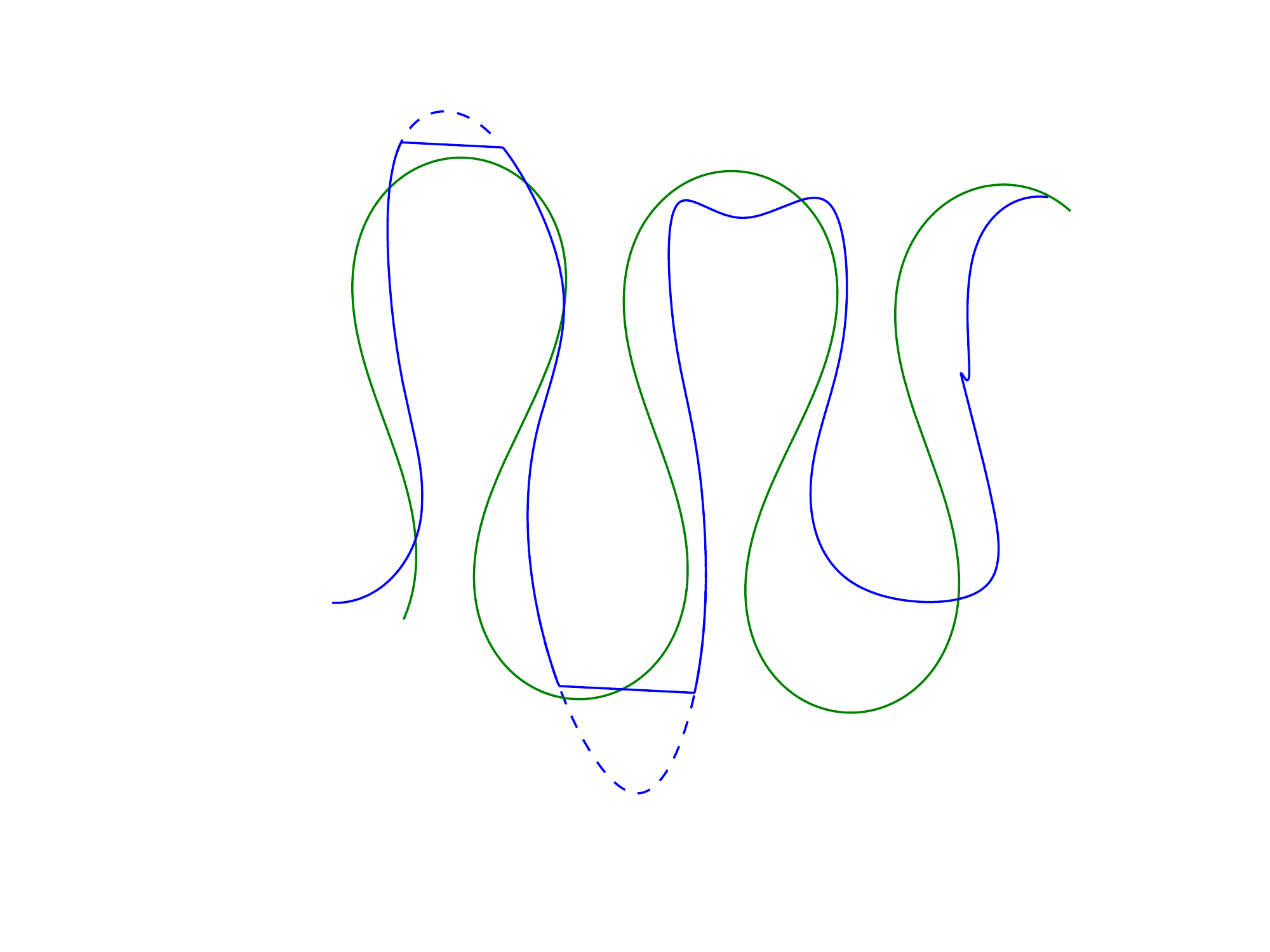}
  \caption{If the given curve moves outside the interval $[u_{\min},u_{\max}]$ (dotted segment), it is simply cut off in these regions.
 The resulting elastica is shown in green. On the left all oscillations of the input curve are counted, on the right the two very small ones are
ignored.}
  \label{fig:ripples}
\end{figure}

At this point we need to take into account the fact that the input curve is not necessarily an elastica.
 For an elastica, we could easily count the oscillations,
but for an arbitrary curve there may be oscillations of different sizes. We find the curve segments where $u$ is monotone,
but we only count such a segment as an oscillation if it has some minimal height:  we have used
 half of the difference $u_{\max}-u_{\min}$ (as defined in Section \ref{sec:elastic-parameters})
 as this minimum.
Moreover, the right hand side of \eqref{eq:cn_umax} need not be between $-1$ and $1$,
or, in the noninflectional case, between $\sqrt{1-1/k^2}$ and $1$.
 We have circumvented this problem by replacing too small values
by $-1$ (or $\sqrt{1-1/k^2}$) and too large values by $1$. The two issues are illustrated in Figure~\ref{fig:ripples}.

We can thus find $s_0$ and $\ell$. We can judge the validity by calculating
\begin{equation*}
  R_3=\frac{1}{L}\int_{u(t)\notin[u_{\min},u_{\max}]}\frac{\ds}{\dt}\,\dt\,.
\end{equation*}

\begin{figure}[here]
\centering
\unitlength=.25\textwidth
\includegraphics[width=\unitlength,clip,trim=20mm 10mm 20mm 10mm]{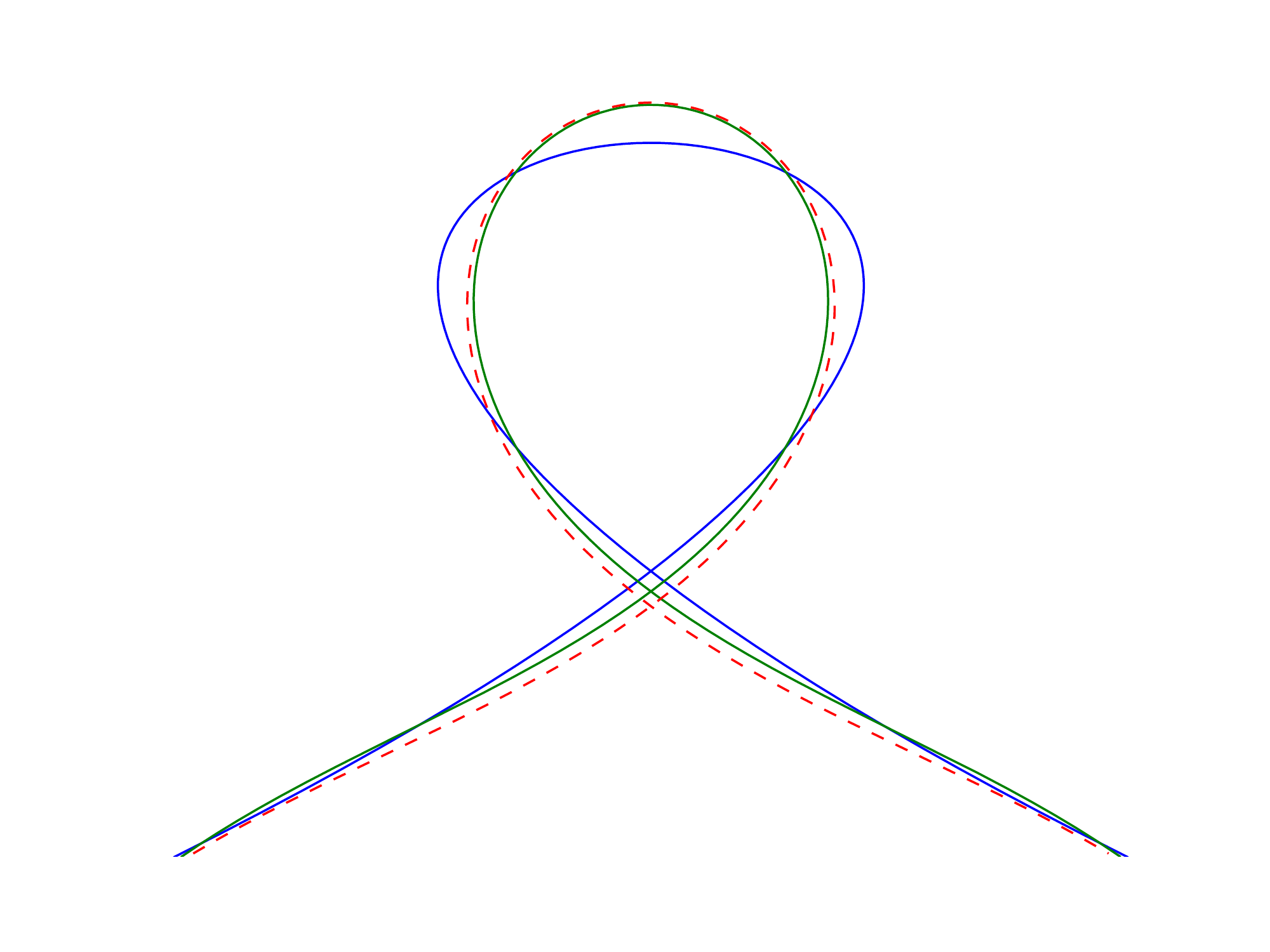}\hfill
\includegraphics[width=\unitlength,clip,trim=20mm 10mm 20mm 10mm]{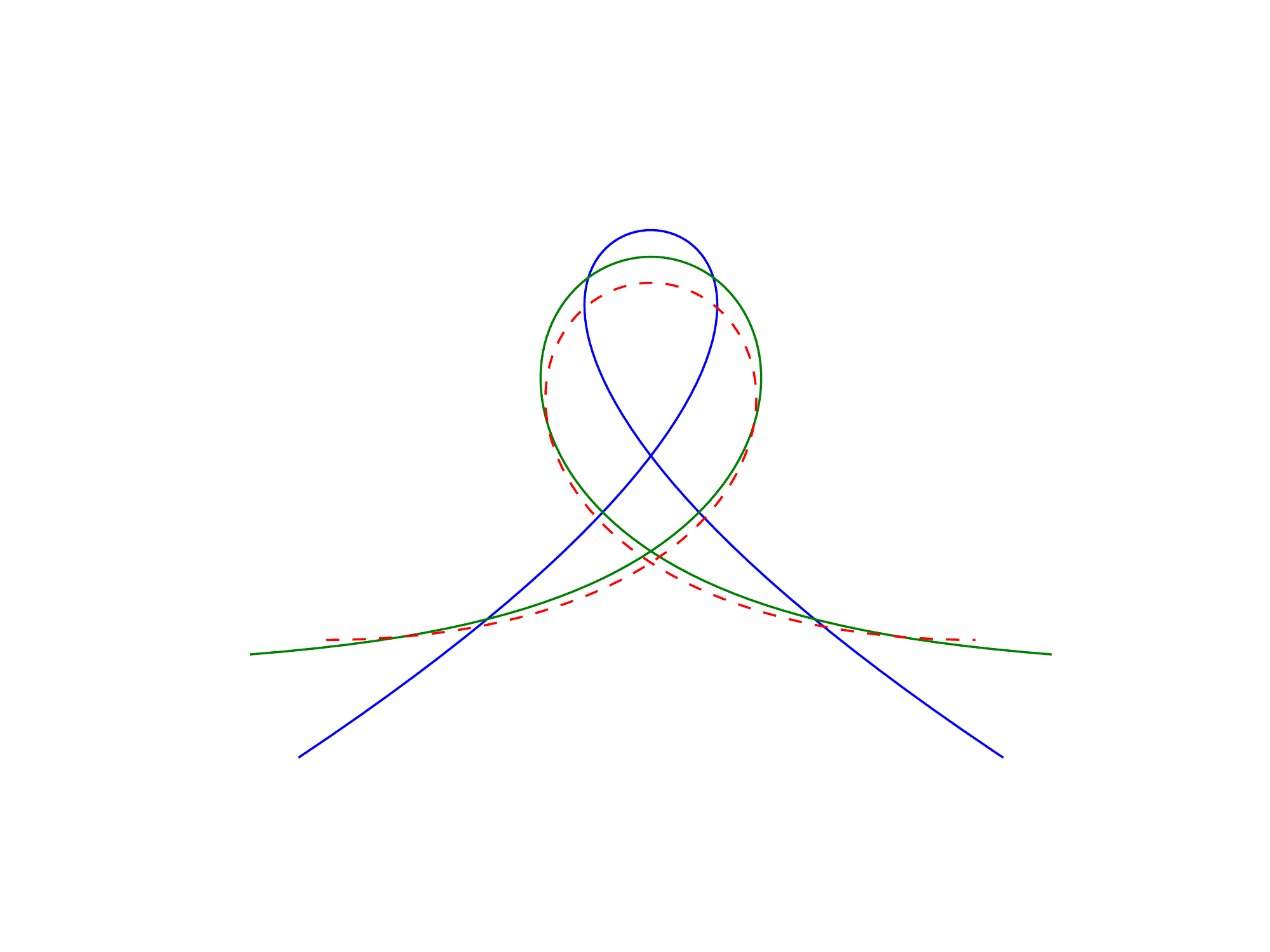}\hfill
\includegraphics[width=\unitlength,clip,trim=20mm 10mm 20mm 10mm]{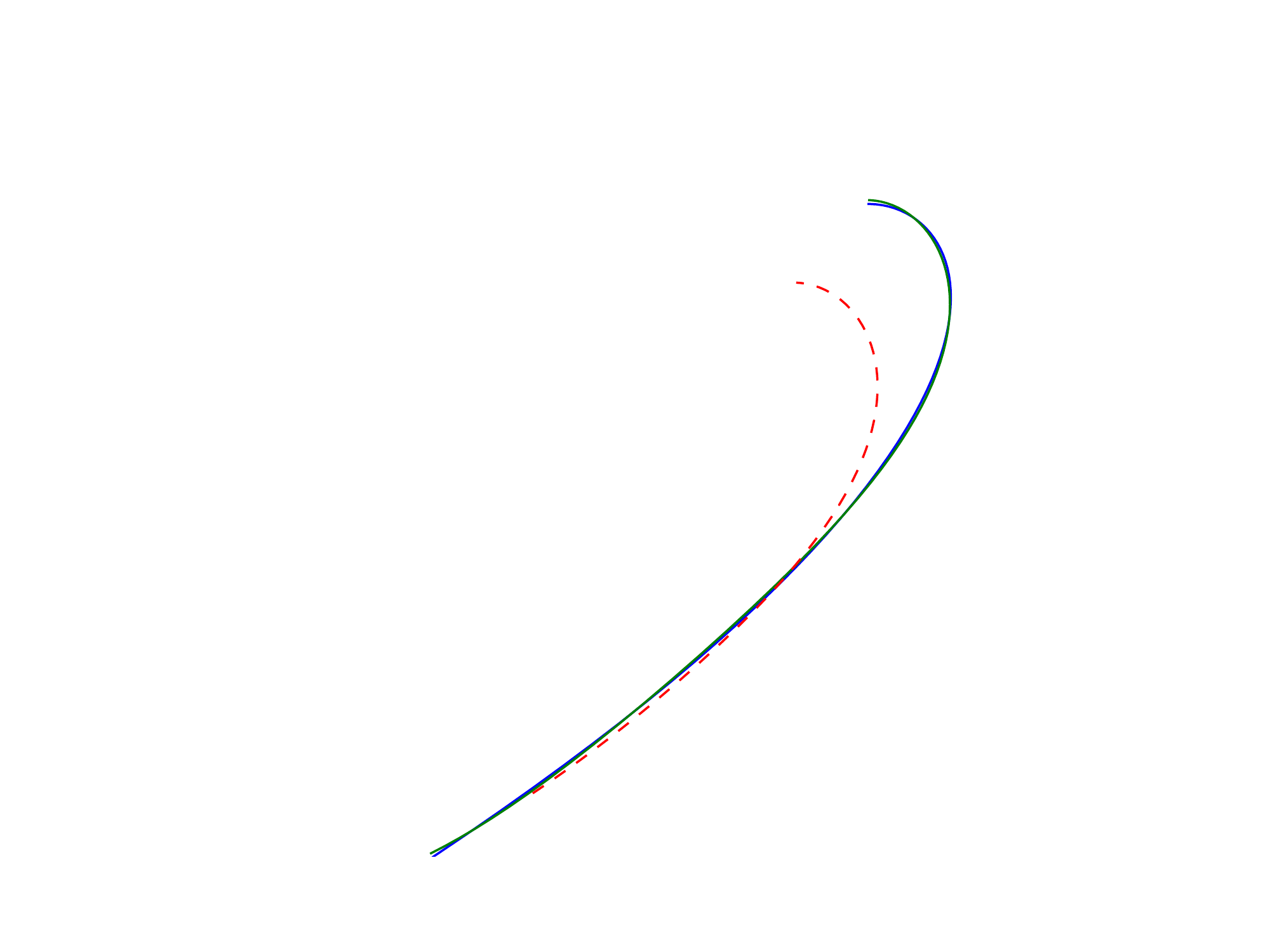} \\
\includegraphics[width=\unitlength,clip,trim=20mm 10mm 20mm 10mm]{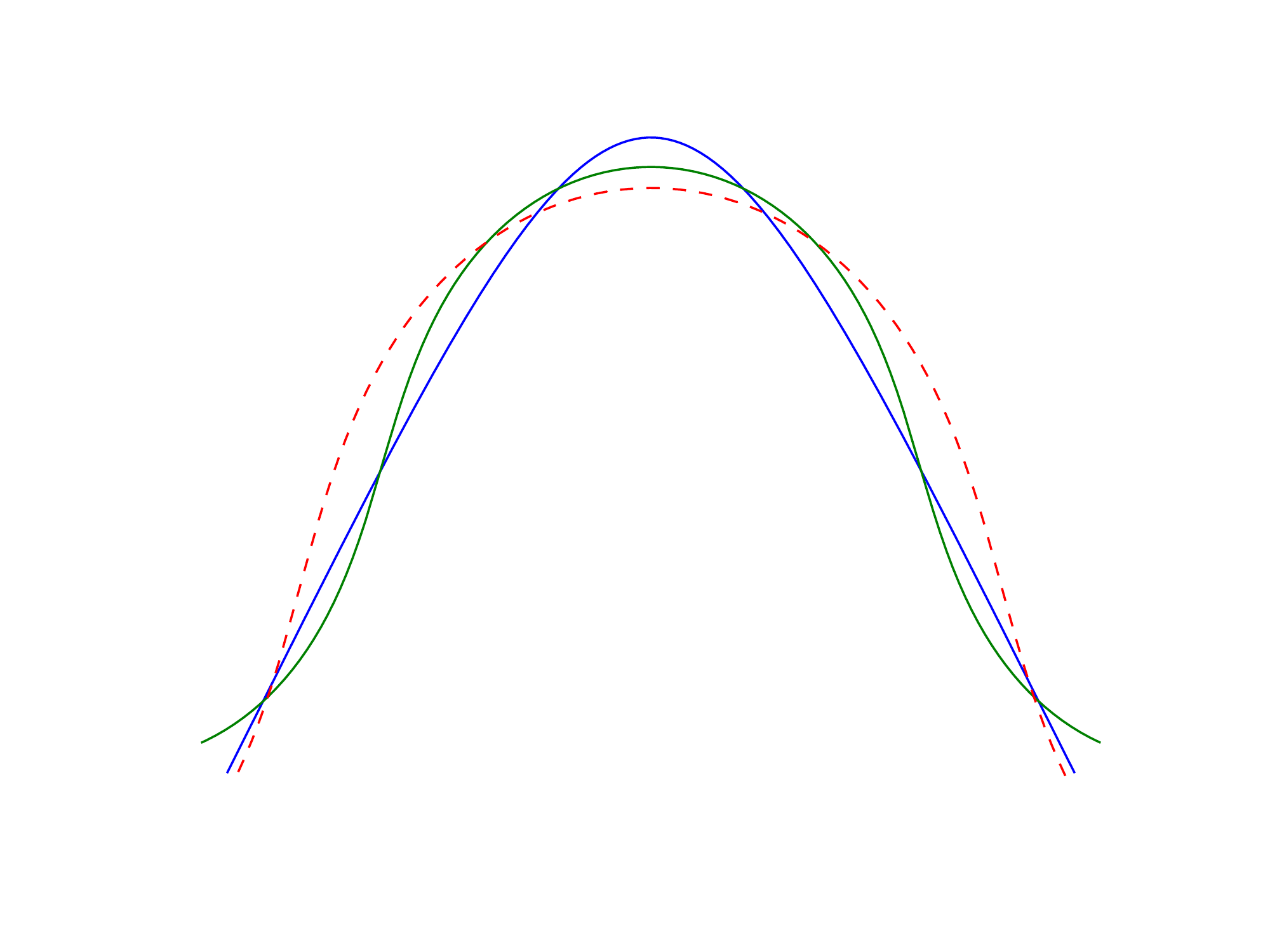}\hfill
\includegraphics[width=\unitlength,clip,trim=20mm 10mm 20mm 10mm]{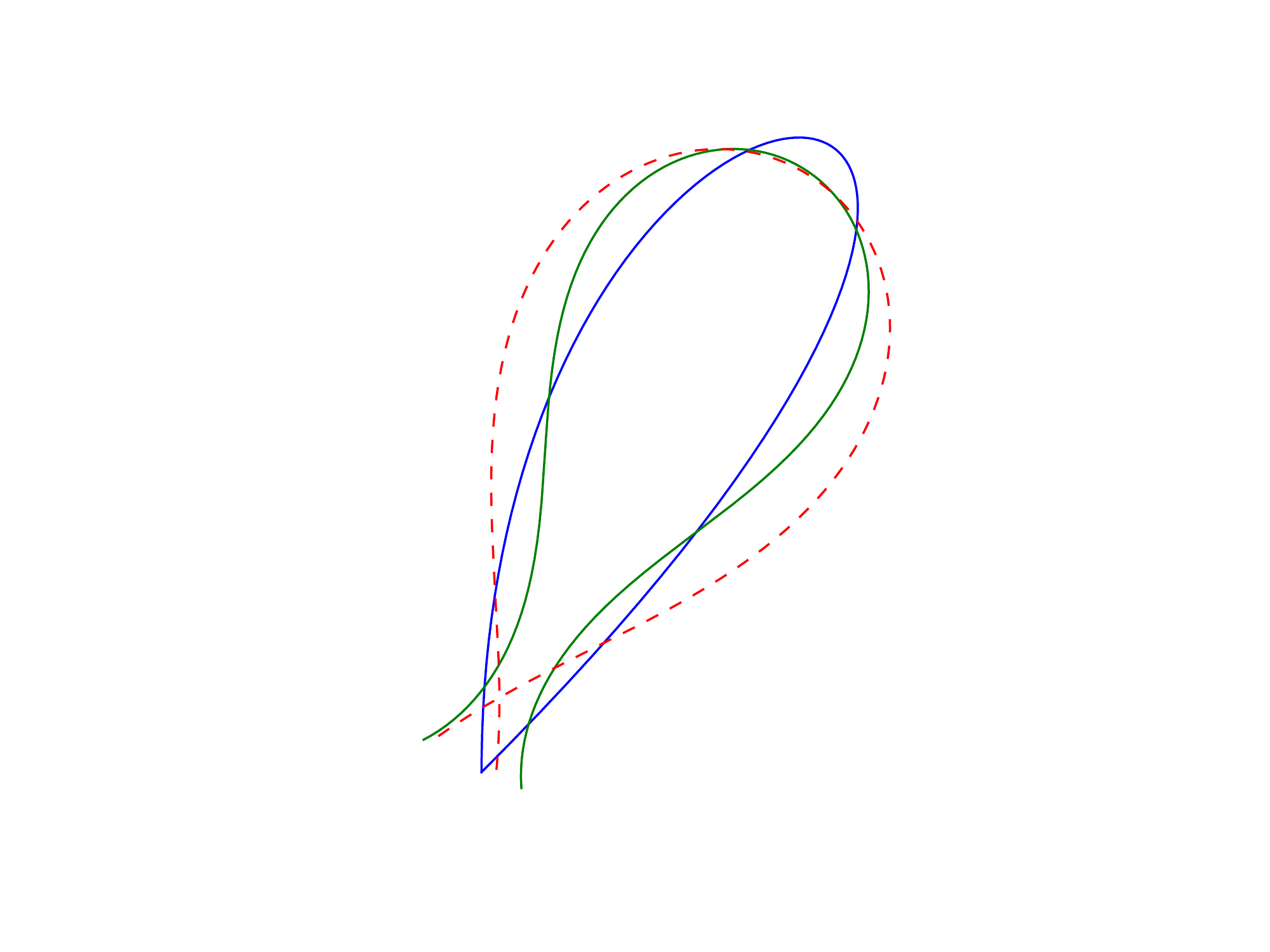}\hfill
\includegraphics[width=\unitlength,clip,trim=20mm 10mm 20mm 10mm]{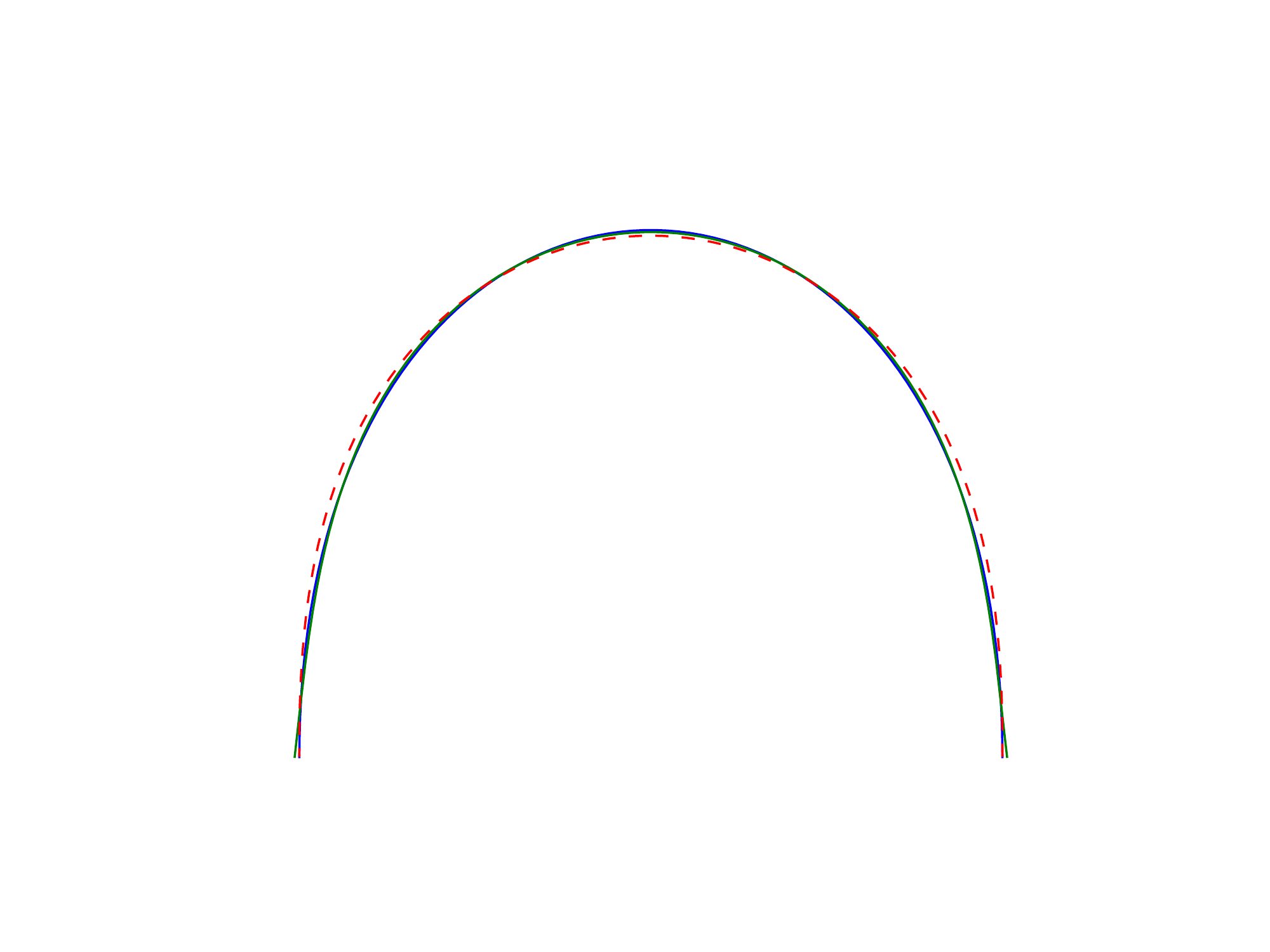} \\
\includegraphics[width=\unitlength,clip,trim=20mm 50mm 20mm 50mm]{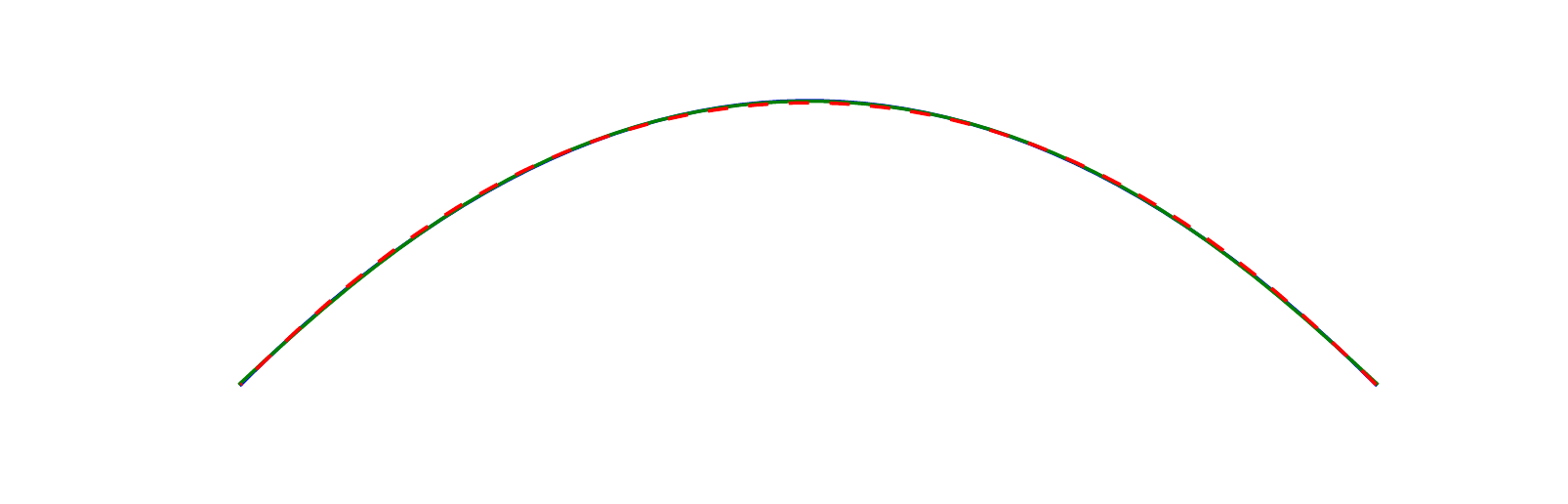}\hfill
\includegraphics[width=\unitlength,clip,trim=20mm 50mm 20mm 50mm]{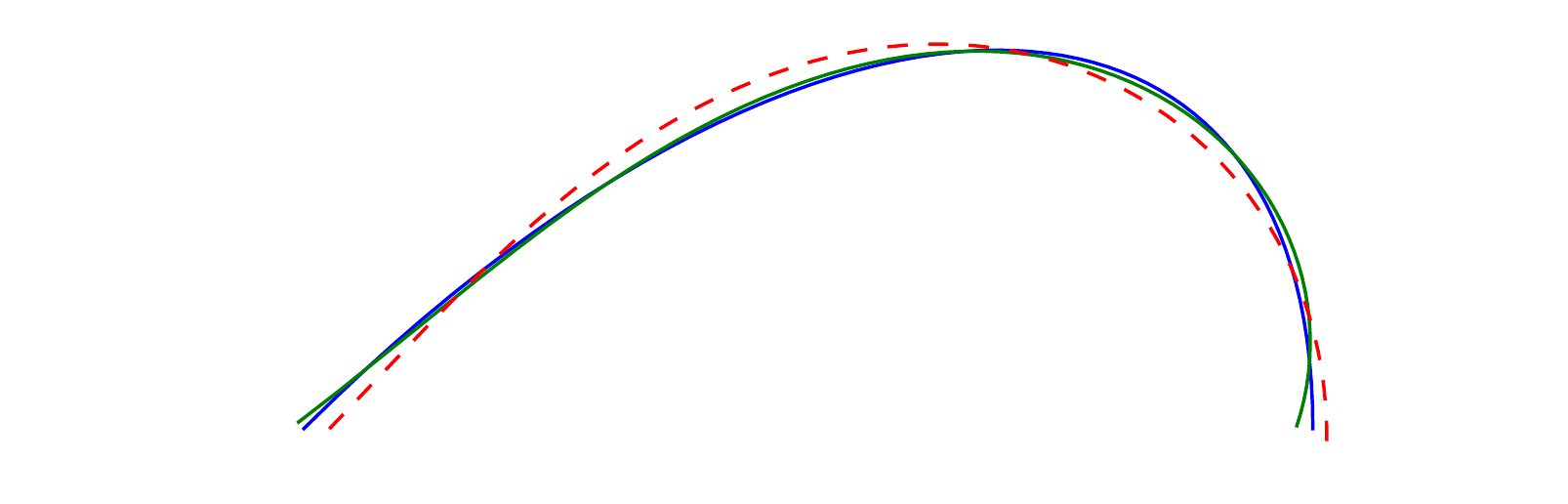}\hfill
\includegraphics[width=\unitlength,clip,trim=20mm 50mm 20mm 50mm]{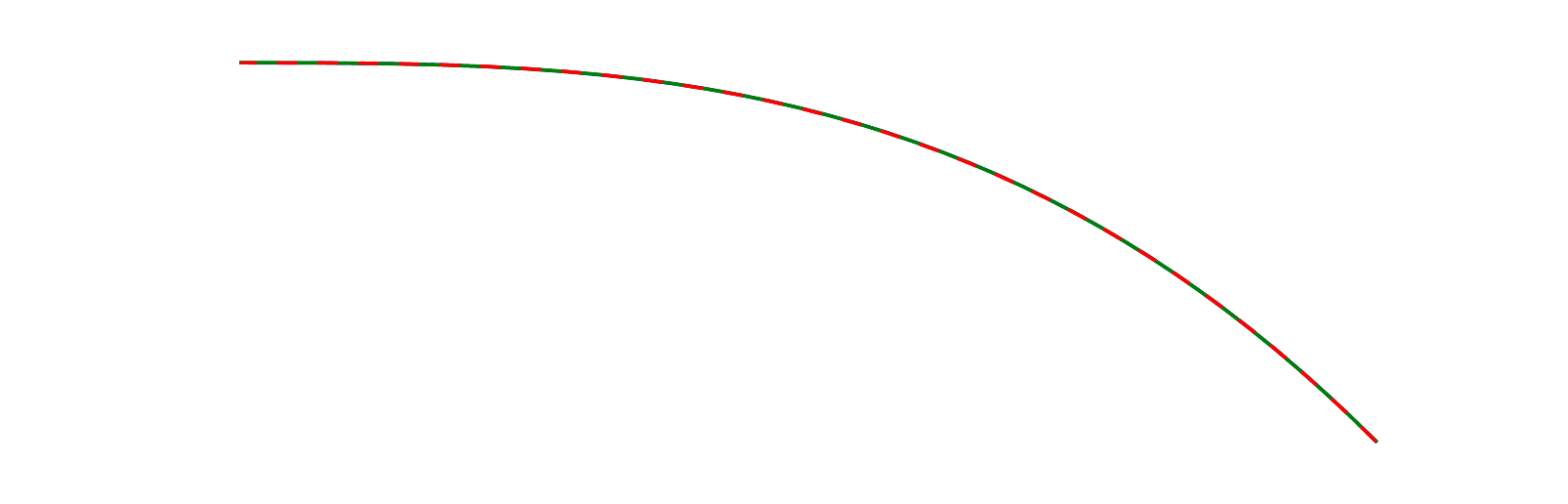} \\
\includegraphics[width=\unitlength,clip,trim=20mm 20mm 20mm 20mm]{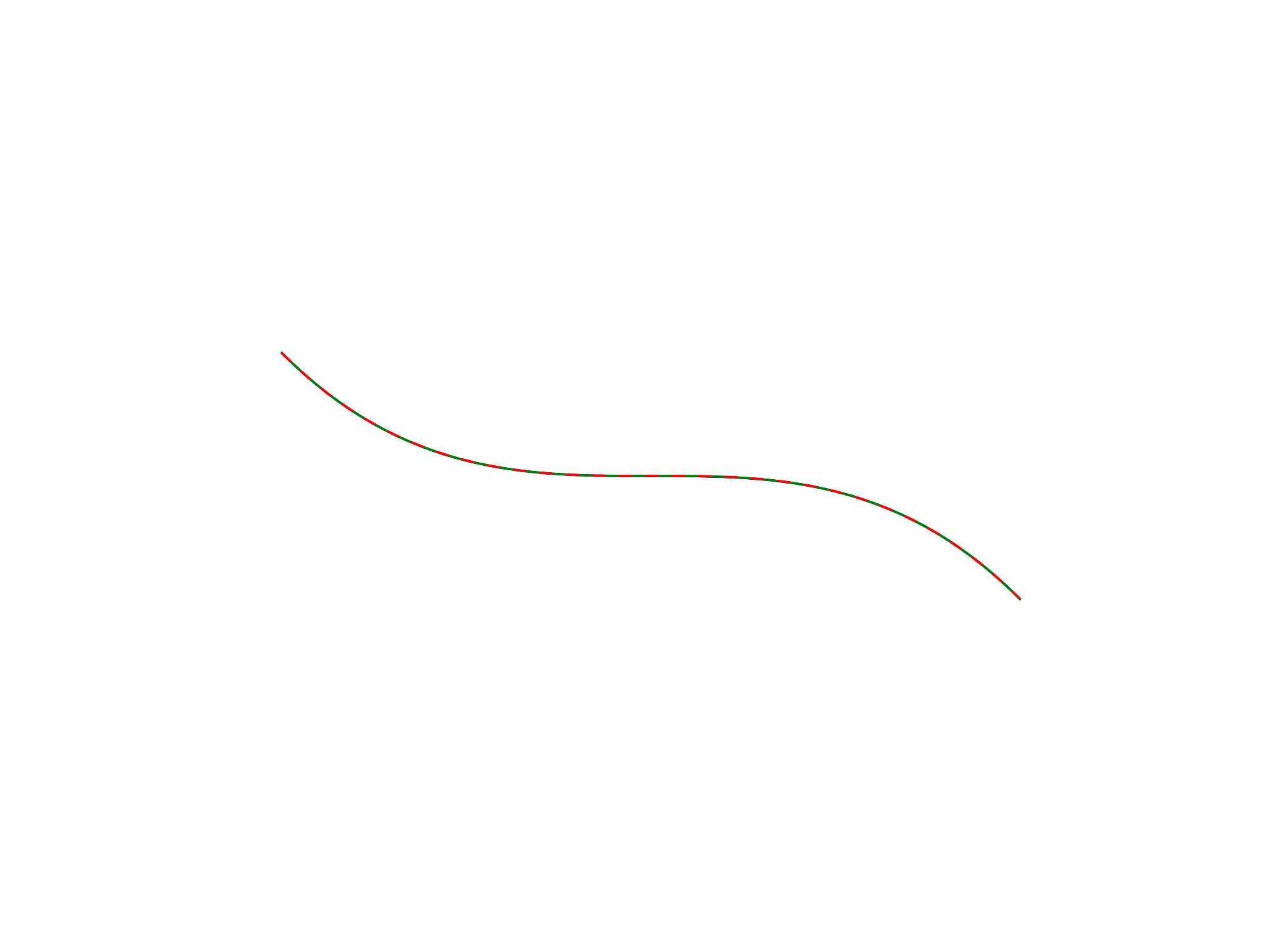}\hfill
\includegraphics[width=\unitlength,clip,trim=20mm 20mm 20mm 20mm]{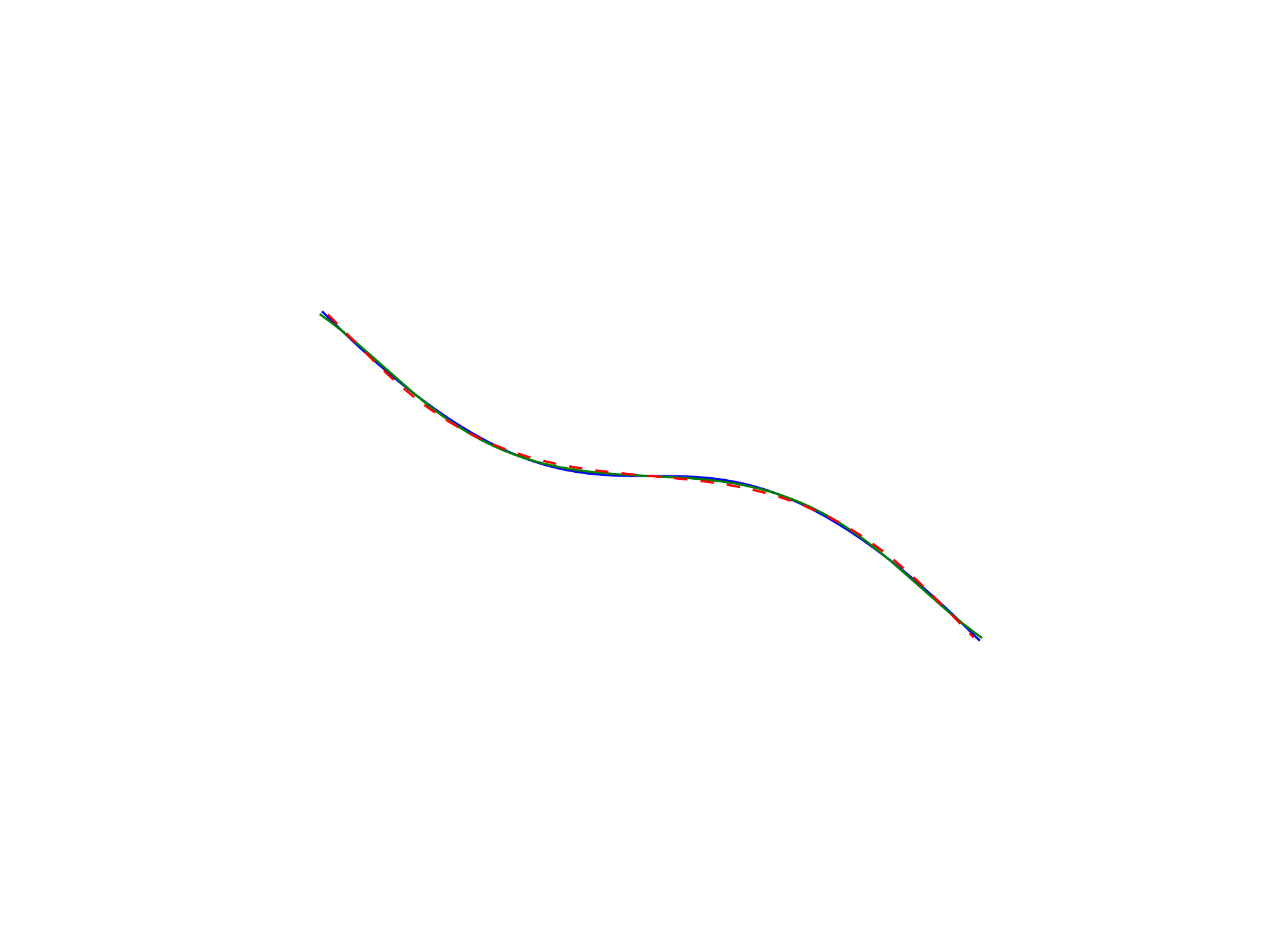}\hfill
\includegraphics[width=\unitlength,clip,trim=20mm 20mm 20mm 20mm]{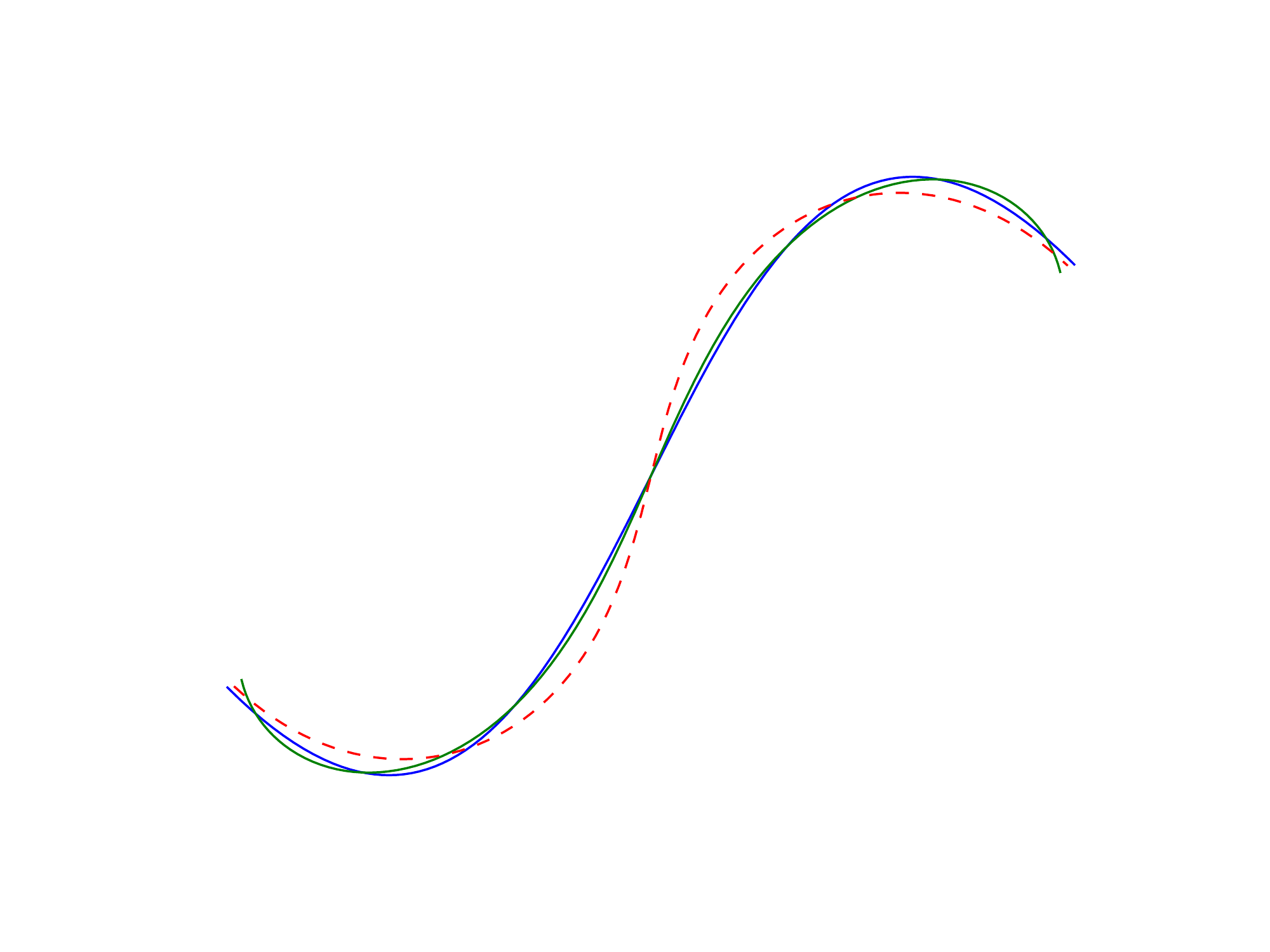}
\caption{Examples of cubic B\'{e}zier curves approximated with elastica. The
solid blue line is the B\'{e}zier curve and the dashed red line is the initial guess for an approximating elastica. The solid green curve is the best approximating elastica found with IPOPT optimization.}
\label{fig:opt1}
\end{figure}

\begin{figure}[here]
\centering
\unitlength=.25\textwidth
\includegraphics[width=\unitlength,clip,trim=20mm 10mm 20mm 10mm]{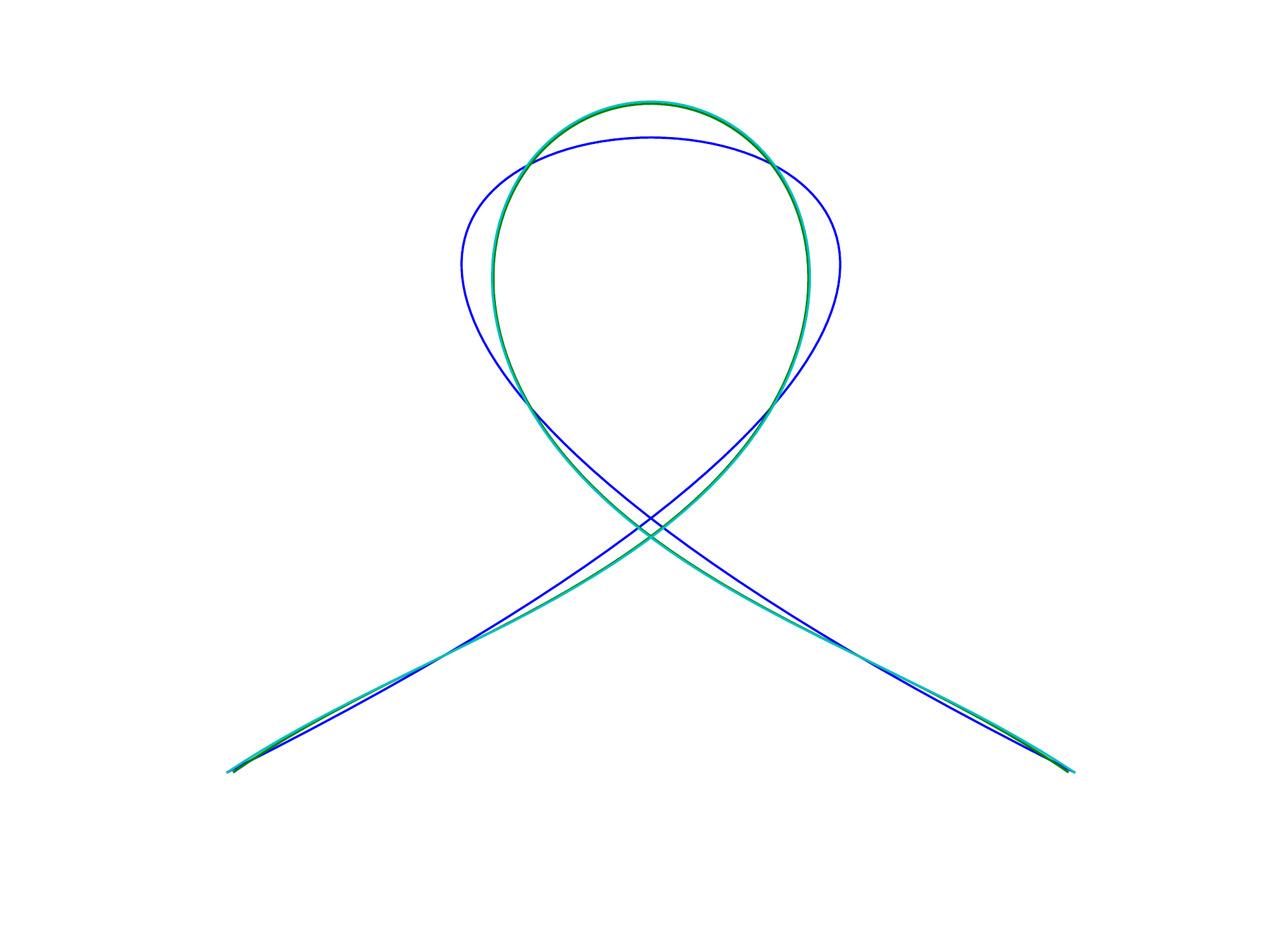}\hfill
\includegraphics[width=\unitlength,clip,trim=20mm 10mm 20mm 10mm]{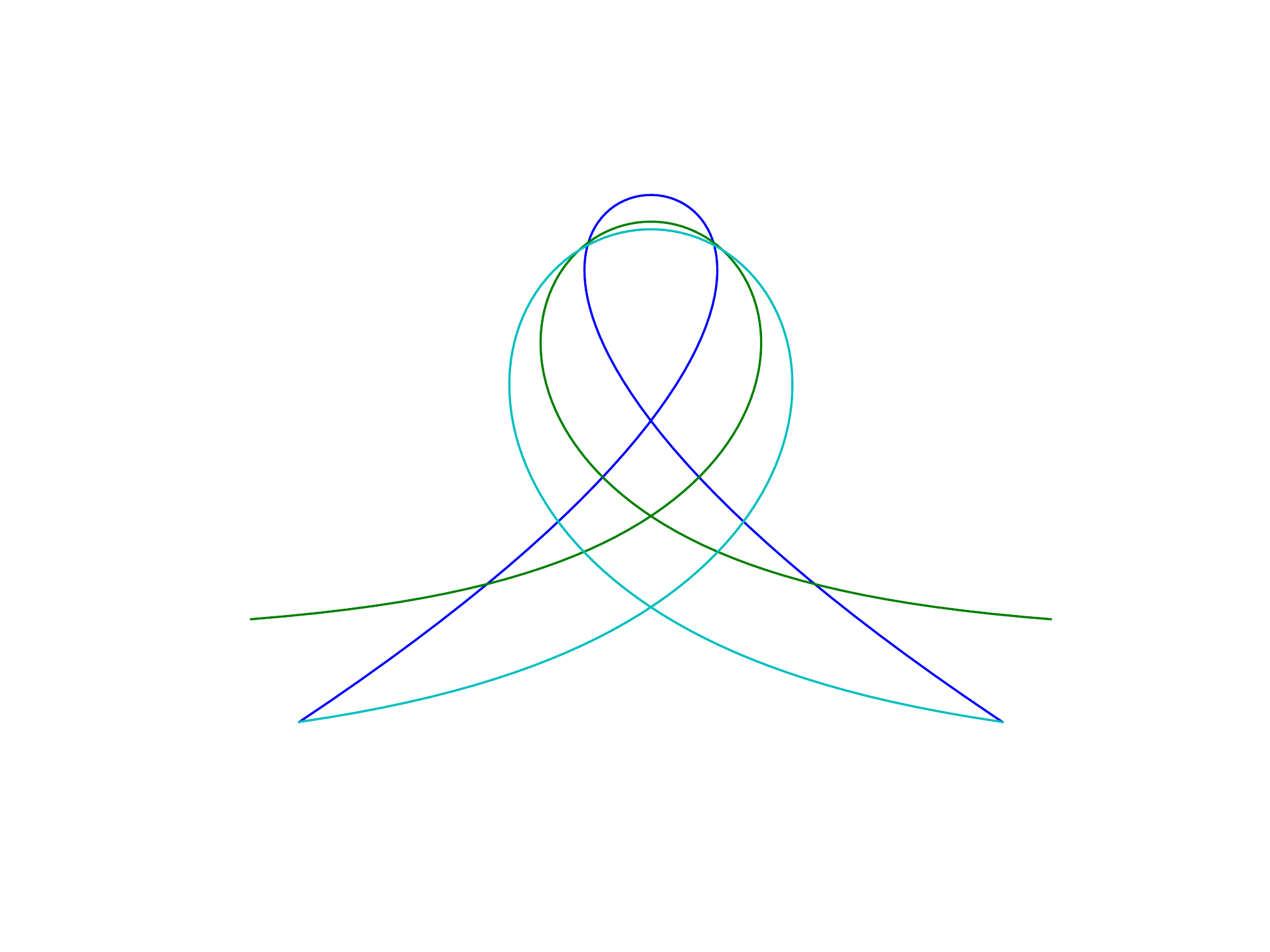}\hfill
\includegraphics[width=\unitlength,clip,trim=20mm 10mm 20mm 10mm]{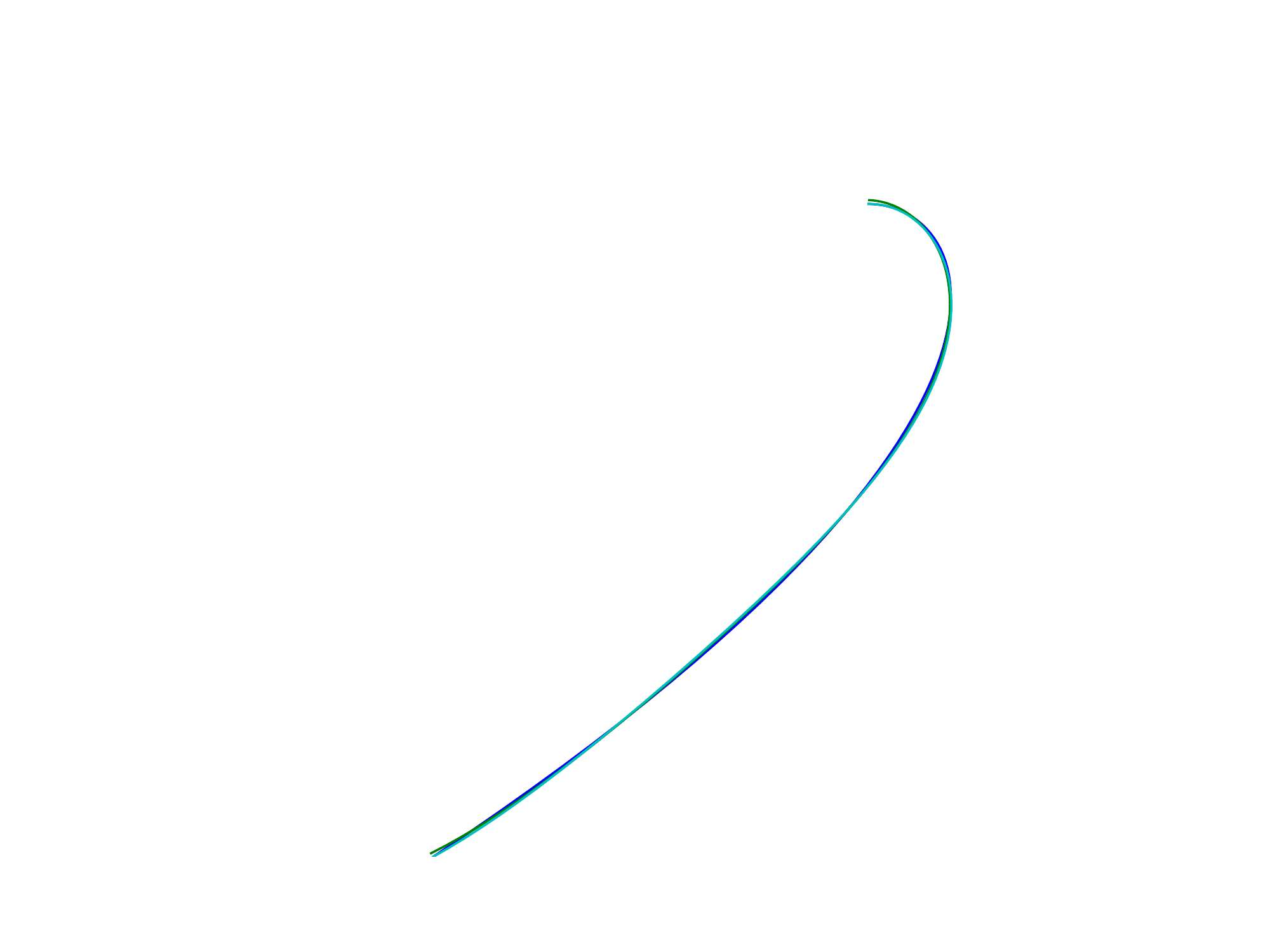} \\
\includegraphics[width=\unitlength,clip,trim=20mm 10mm 20mm 10mm]{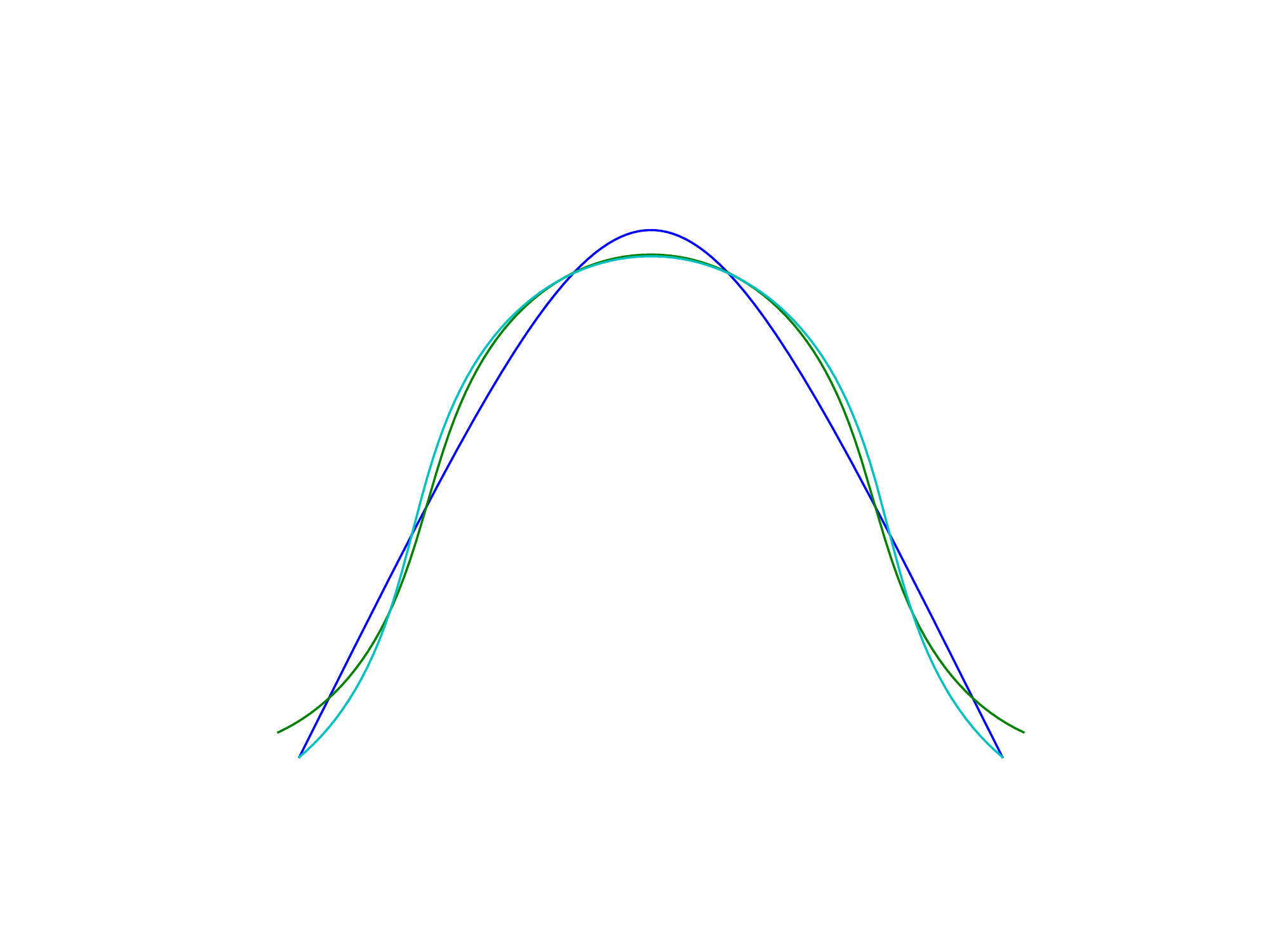}\hfill
\includegraphics[width=\unitlength,clip,trim=20mm 10mm 20mm 10mm]{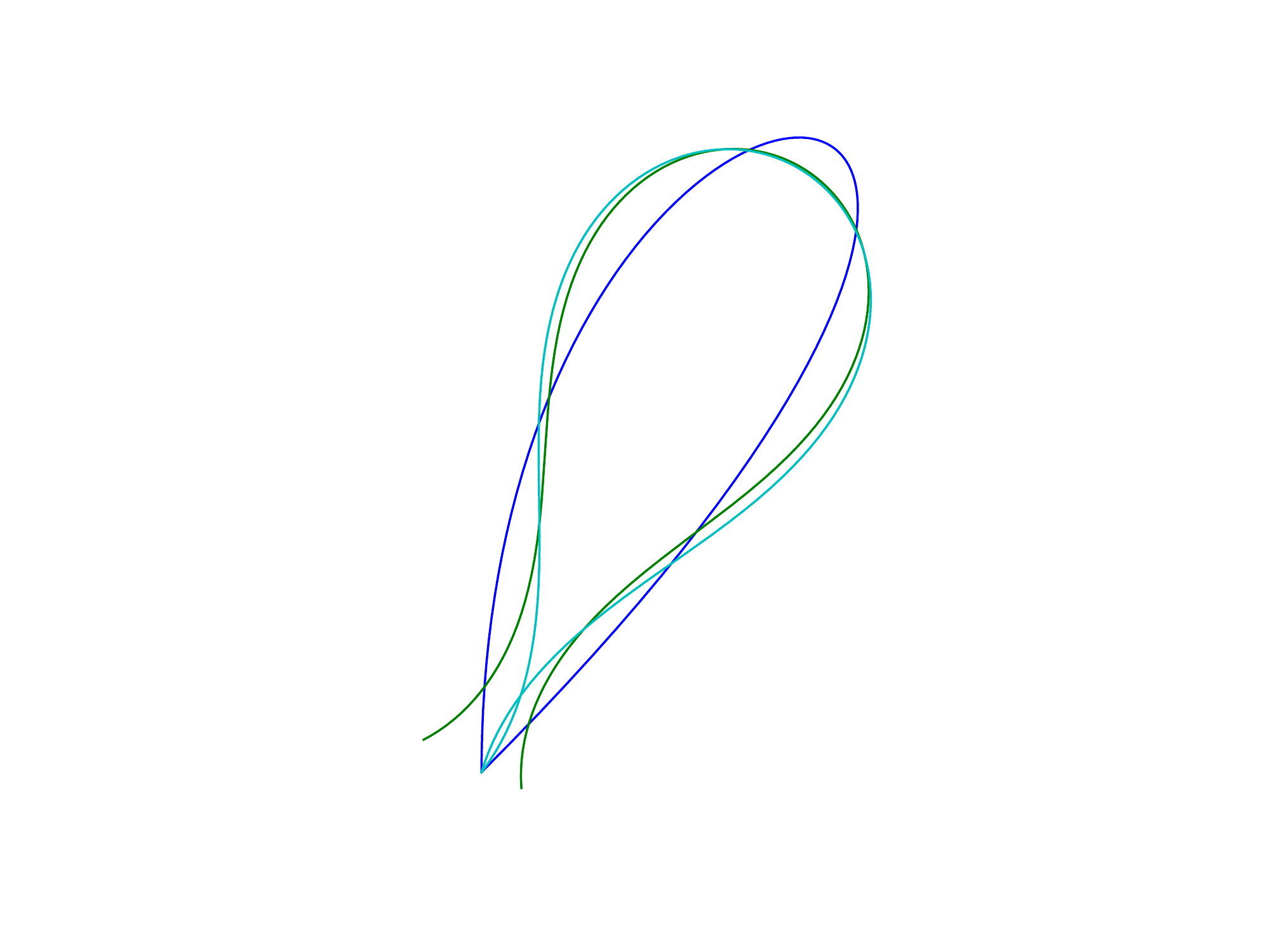}\hfill
\includegraphics[width=\unitlength,clip,trim=20mm 10mm 20mm 10mm]{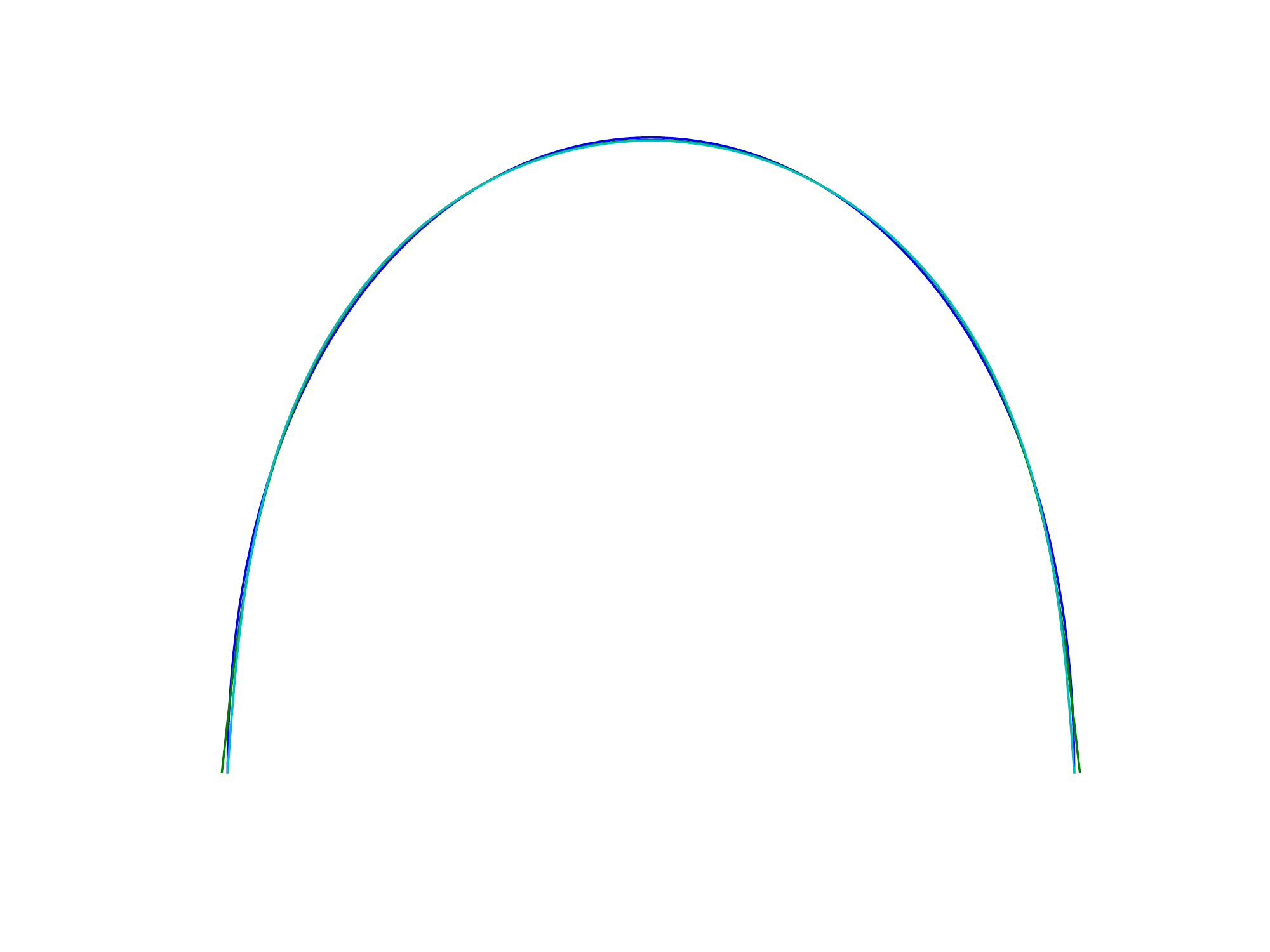} \\
\includegraphics[width=\unitlength,clip,trim=20mm 50mm 20mm 50mm]{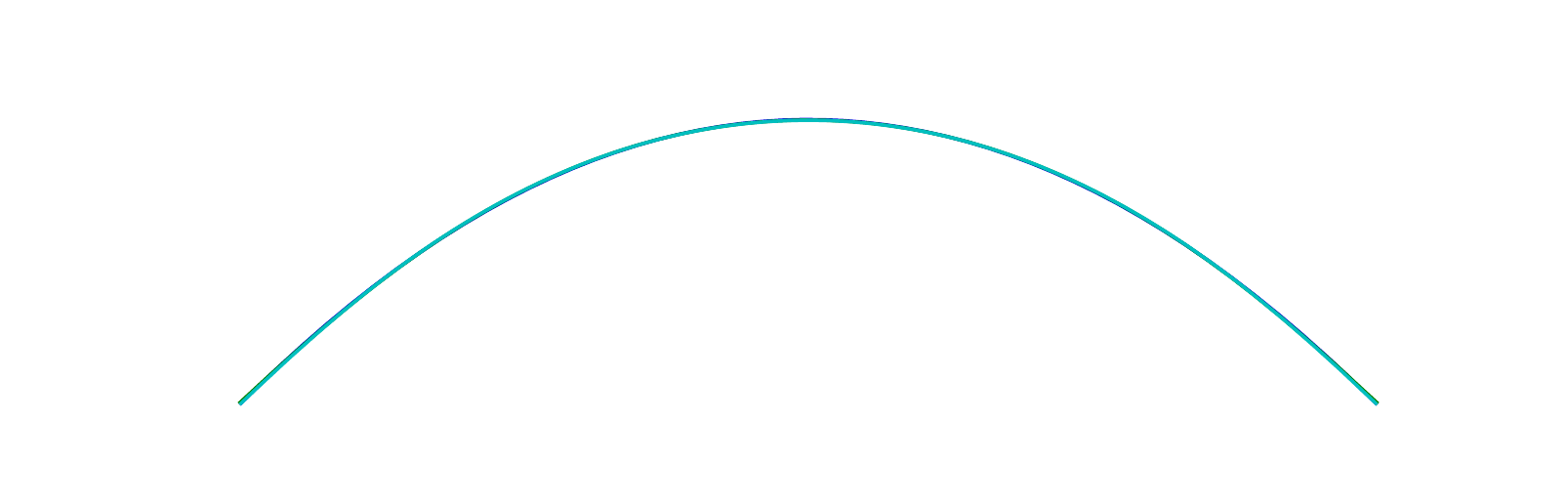}\hfill
\includegraphics[width=\unitlength,clip,trim=20mm 50mm 20mm 50mm]{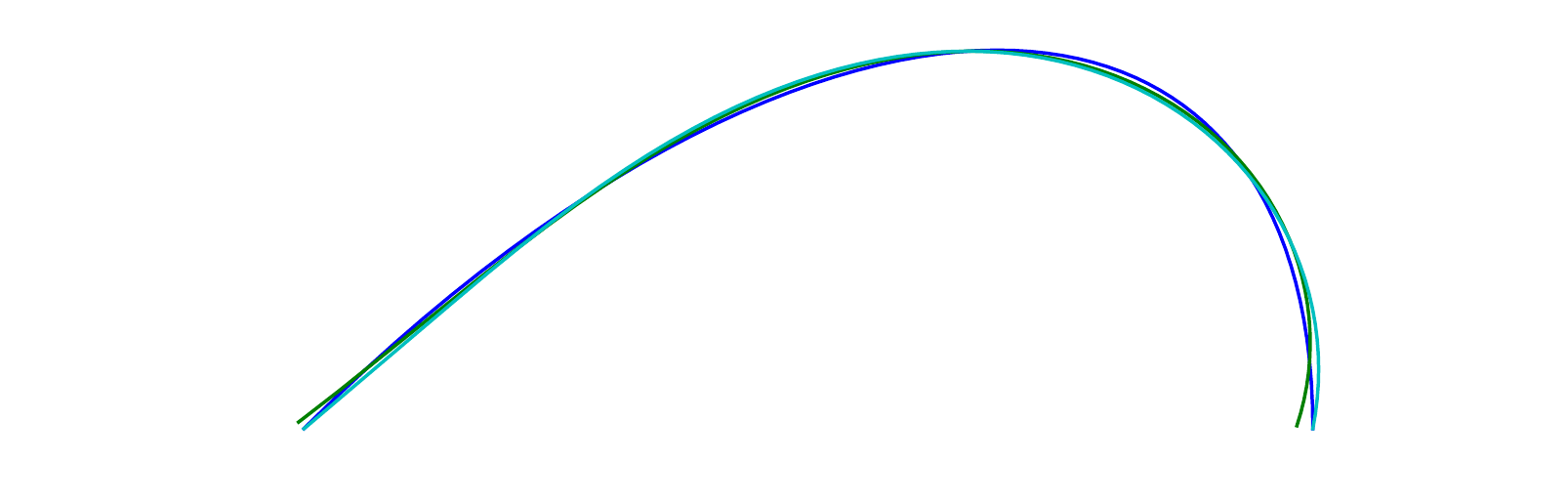}\hfill
\includegraphics[width=\unitlength,clip,trim=20mm 50mm 20mm 50mm]{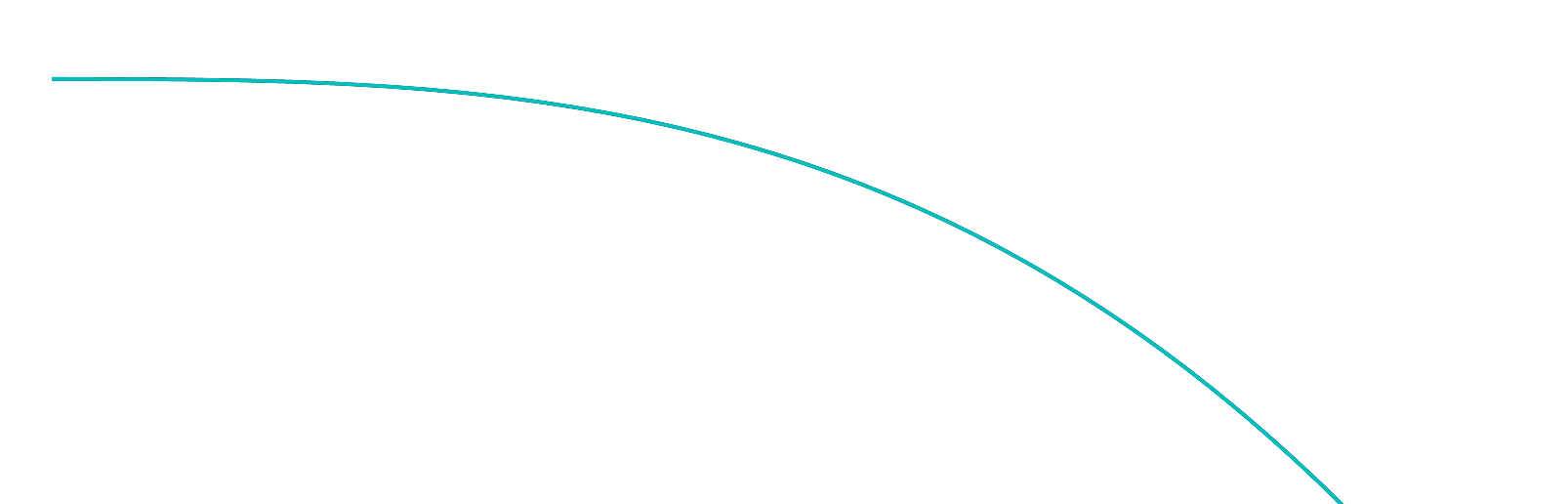} \\
\includegraphics[width=\unitlength,clip,trim=20mm 20mm 20mm 20mm]{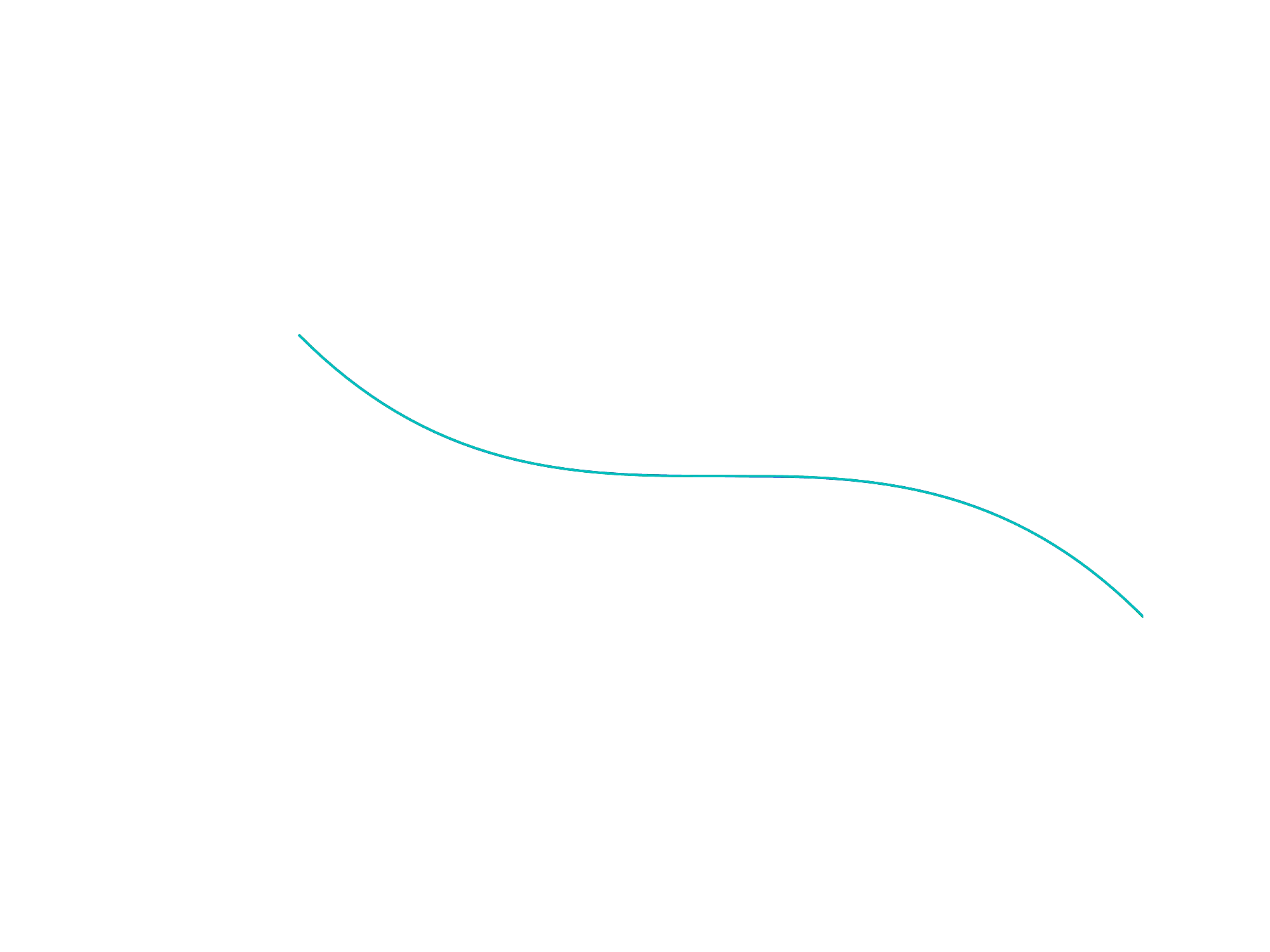}\hfill
\includegraphics[width=\unitlength,clip,trim=20mm 20mm 20mm 20mm]{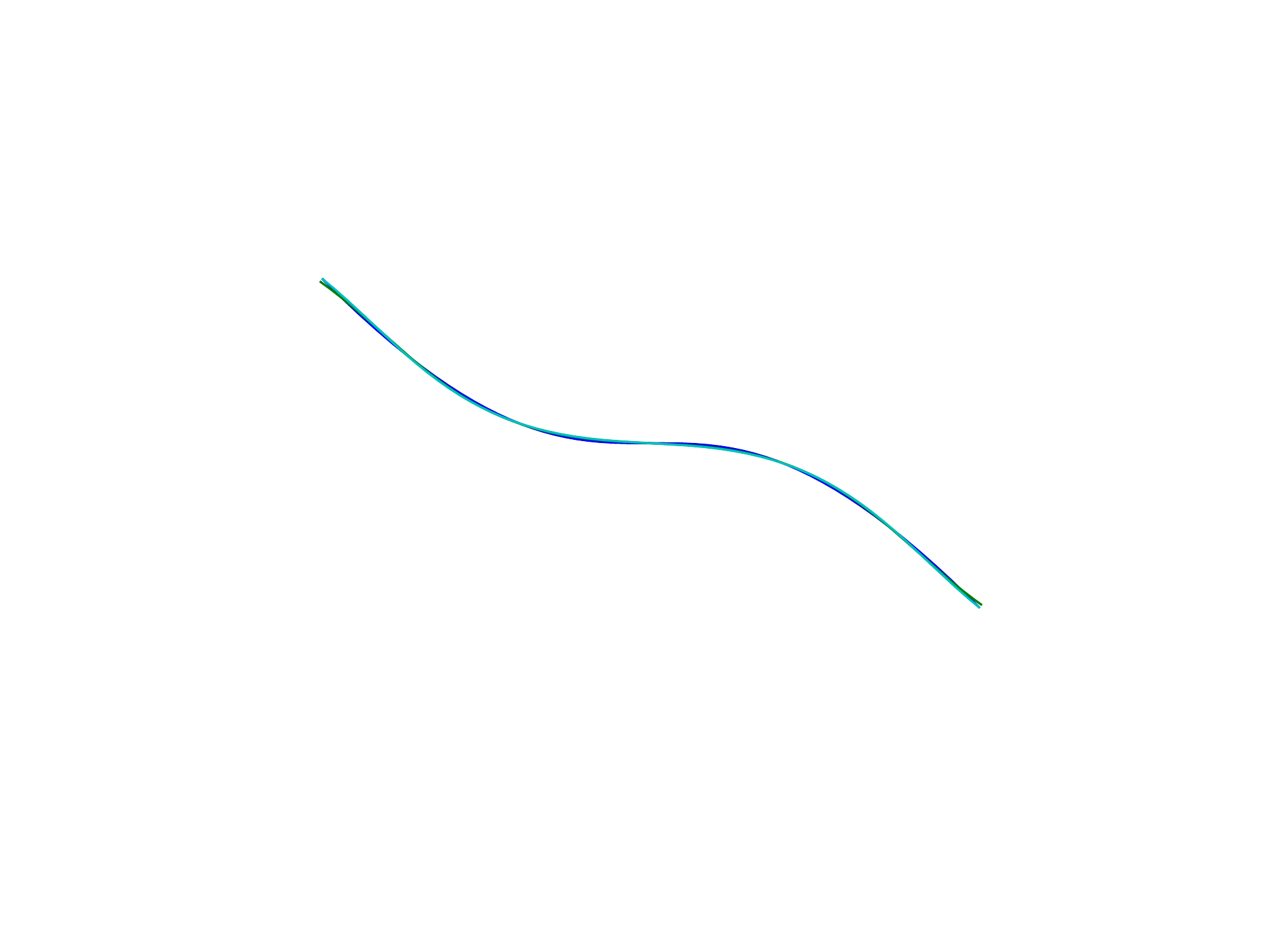}\hfill
\includegraphics[width=\unitlength,clip,trim=20mm 20mm 20mm 20mm]{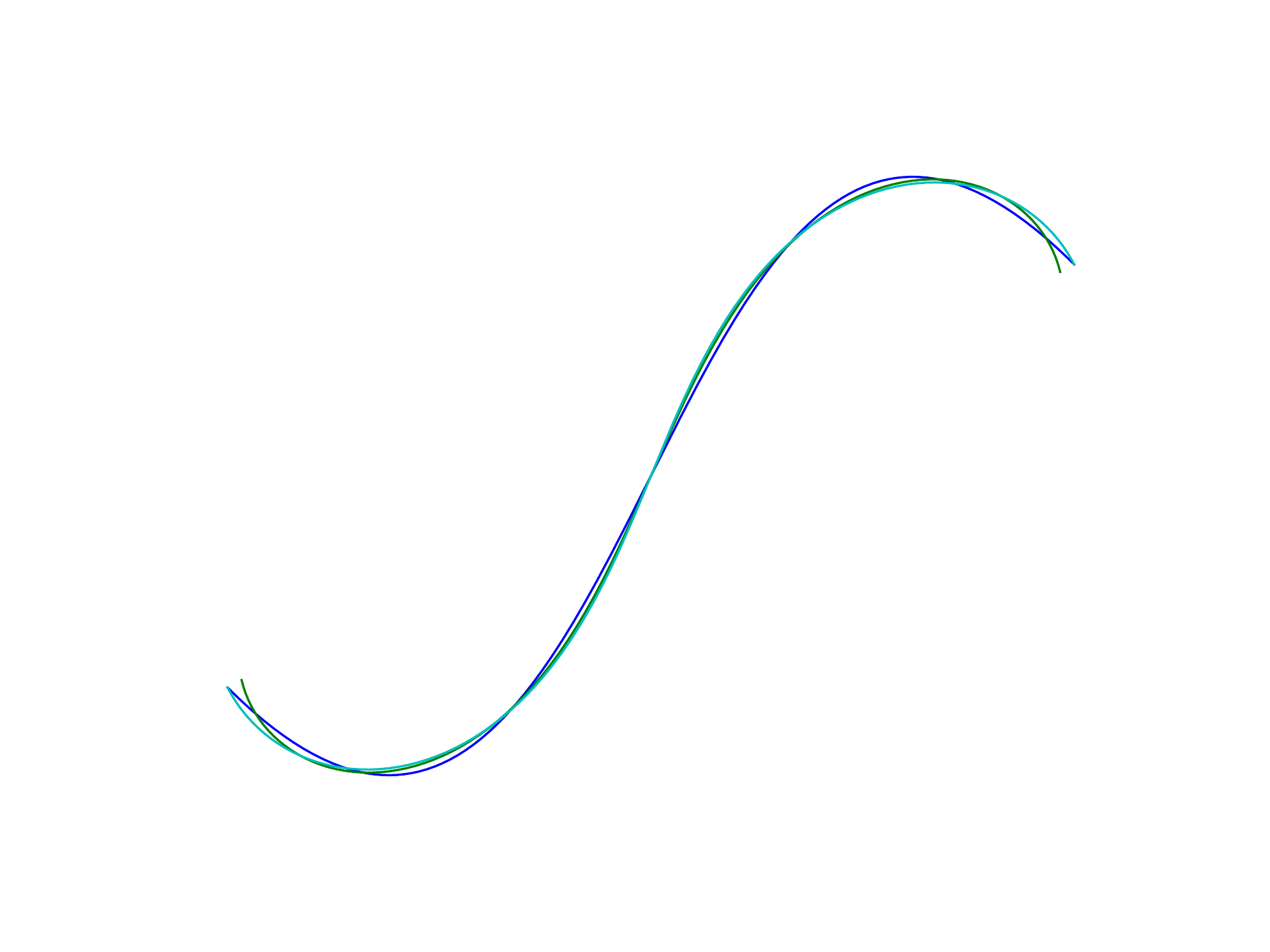}
\caption{The solid cyan curve is the best approximating elastica having the same endpoints as the original curve (blue). The green curve is
the output curve with free endpoints, as in Figure~\ref{fig:opt1}.}
\label{fig:opt2}
\end{figure}

We finally determine the translation by \eqref{eq:translate} and
define the residual as
\begin{equation*}
  R_4(\mathbf{p}_0)=\sqrt{\frac{2}{L^3}\mathcal{F}(\mathbf{p}_0)}\,,
\end{equation*}
where $\mathbf{p}_0$ is the vector of control parameters found by the above procedure.

We have tested the procedure on a selection of cubic B\'{e}zier curves, displayed in Figures~\ref{fig:opt1} and \ref{fig:opt2}. In Table~\ref{tab:Bezier} we have reported the residuals.

\begin{table}[here]
  \caption{The first column refers to the examples in Figure~\ref{fig:opt1}, the next three report the residuals
    $R_1$, $R_2$, and $R_3$, in the approximation process.
    Next, we have the normalized $L^2$-distance, $R_4=\sqrt{2\mathcal{F}/L^3}$, for the initial
    guess and then $R_4$, the gradient norm $\|\nabla\mathcal{F}\|$ and the number of iterations for the optimized elastica without and with endpoint constraints. We use the same initial guess, whether we constrain the endpoints or not.}
  \label{tab:Bezier}
  \centering
	\begin{small}
  \begin{tabular}{cccccccccc}
    \hline\noalign{\smallskip}
    & $R_1$ & $R_2$ &
    $R_3$
    & $R_4(\mathbf{p}_0)$ & $R_4(\mathbf{p}_{\mathrm{opt}})$ & $\|\nabla\mathcal{F}(\mathbf{p}_{\mathrm{opt}})\|$ & $\sharp$ iter & $R_4(\mathbf{p}_{\mathrm{opt}}^*)$ & $\sharp$ iter$^*$
    \\
		\noalign{\smallskip}\hline\noalign{\smallskip}
1 & 0.46 & 0.14 & 0.0 & 0.0097 & 0.0080 & 2.2e-08 & 35 & 0.0081 & 10 \\
2 & 0.45 & 0.23 & 0.52 & 0.048 & 0.038 & 5.1 & 1000\dag & 0.063 & 1000\dag \\
3 & 0.14 & 0.020 & 0.19 & 0.077 & 0.0018 & 6.0e-09 & 26 & 0.0023 & 16 \\
4 & 0.68 & 0.17 & 0.14 & 0.022 & 0.012 & 5.9e-09 & 8 & 0.014 & 9 \\
5 & 0.65 & 0.27 & 0.14 & 0.031 & 0.018 & 5.2e-09 & 9 & 0.021 & 8 \\
6 & 0.099 & 0.025 & 0.048 & 0.0036 & 0.0011 & 8.1e-11 & 24 & 0.0015 & 14 \\
7 & 0.064 & 0.0044 & 0.048 & 0.0010 & 0.00032 & 1.7e-10 & 20 & 0.00041 & 104 \\
8 & 0.27 & 0.069 & 0.14 & 0.011 & 0.0032 & 1.5e-09 & 199 & 0.0046 & 100 \\
9 & 0.017 & 0.00033 & 0.0 & 0.00013 & 0.0012 & 0.053 & 1000 & 5.1e-05 & 15 \\
10 & 0.020 & 0.00062 & 0.0 & 0.00012 & 9.9e-05 & 1.2e-09 & 165 & 0.00011 & 29 \\
11 & 0.33 & 0.014 & 0.0 & 0.0050 & 0.0017 & 8.3e-11 & 76 & 0.0020 & 11 \\
12 & 0.36 & 0.10 & 0.19 & 0.015 & 0.0041 & 7.3e-09 & 178 & 0.0053 & 83 \\
 \noalign{\smallskip}\hline
  \end{tabular}
	\end{small}
  \flushleft
  \dag IPOPT terminated because iteration count reached maximum (which was set to 1000). 
\end{table}


\section {Conclusions, discussion and  related work}  \label{sec:conclusion}
\begin{figure}[here]
\centering
\includegraphics[height=30mm]{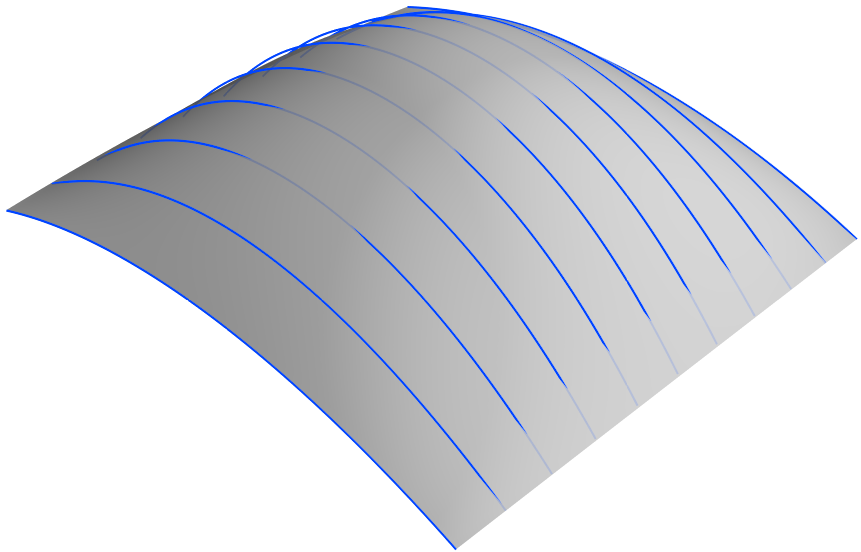} \quad \quad
\includegraphics[height=30mm]{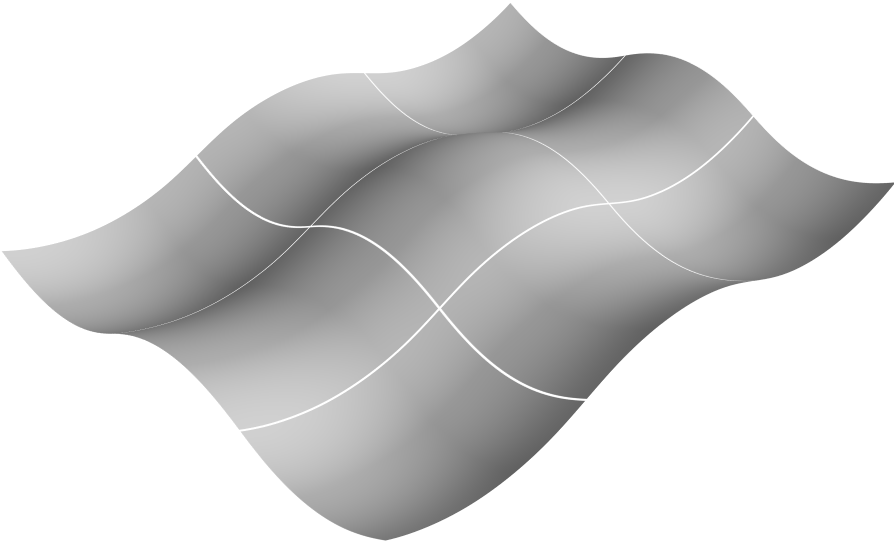}
\caption{Surface rationalization for hot-blade cutting. Left: a CAD surface patch is approximated by a family of
elastic curves.  Right: a surface segmented into elastica-foliated pieces. }
\label{fig:surface_example}
\end{figure}

When combined with a suitable segmentation, we have found that
the method described here gives an effective algorithm for approximating curves by piecewise elastic curves.
The degrees of freedom allow
piecewise elastica, with $C^1$ continuity at the joins if end-points are fixed, and with $C^2$ continuity
if end-points are allowed to move.
The method is incorporated in on-going work on approximating \emph{surfaces} by segmented surfaces, with
 the segments swept out by elastic curves (Figure  \ref{fig:surface_example}), for the purpose of
manufacturing architectural formwork by robotic hot-blade cutting \cite{madifa}, \cite{robarch}.

  Our choice of parameters allows us to work with
  analytic expressions of elastic curves. The advantage is that when
  the seven parameters are known any subsequent calculation is
  accurate and easy to perform. The disadvantage is that the geometric
  meaning is not obvious for all seven parameters. A more geometric
  set of parameters is the length and the end points and tangents; the
  problem with that alternative is that we would have to solve a nonlinear boundary
  problem to determine the curve and, more severely, the solution is not
  unique.  For example, if we imagine that we rotate one of the tangents through
  $2\pi$ and follow a continuous family of elastic curves then the
  curve at the end will be different from the curve at the start.

 Note that all the examples in Section \ref{sec:approximating} have, for convenience, been computed using the
$L^2$ distance. Given that our initial guess method is based on the
curvature,  an $H^2$ norm may be a more natural choice, but in practice more complicated.

  Another approach to the problem is to minimize
  \begin{math}
   \int\|\gamma-x\|\,\ds+\beta\int\kappa^2\,\ds\,.
  \end{math}
  Here it is not clear how to choose $\beta$. If $\beta=0$ then we
  obtain $x=\gamma$ and in the limit where $\beta\to\infty$ we obtain
  the best line segment approximation to $\gamma$. In order to make a
  purely numerical approach work we need to have a good approximation
  to an elastic curve at all times while we minimize the distance to
  the target curve $\gamma$. This means we will have to solve a
  nonlinear equation at each step of the optimization and also
  find the sensitivities of this solution.

\begin{figure}[h]
\centering
\unitlength=0.33\textwidth
\includegraphics[width=\unitlength,trim=35mm 15mm 20mm 25mm,clip]{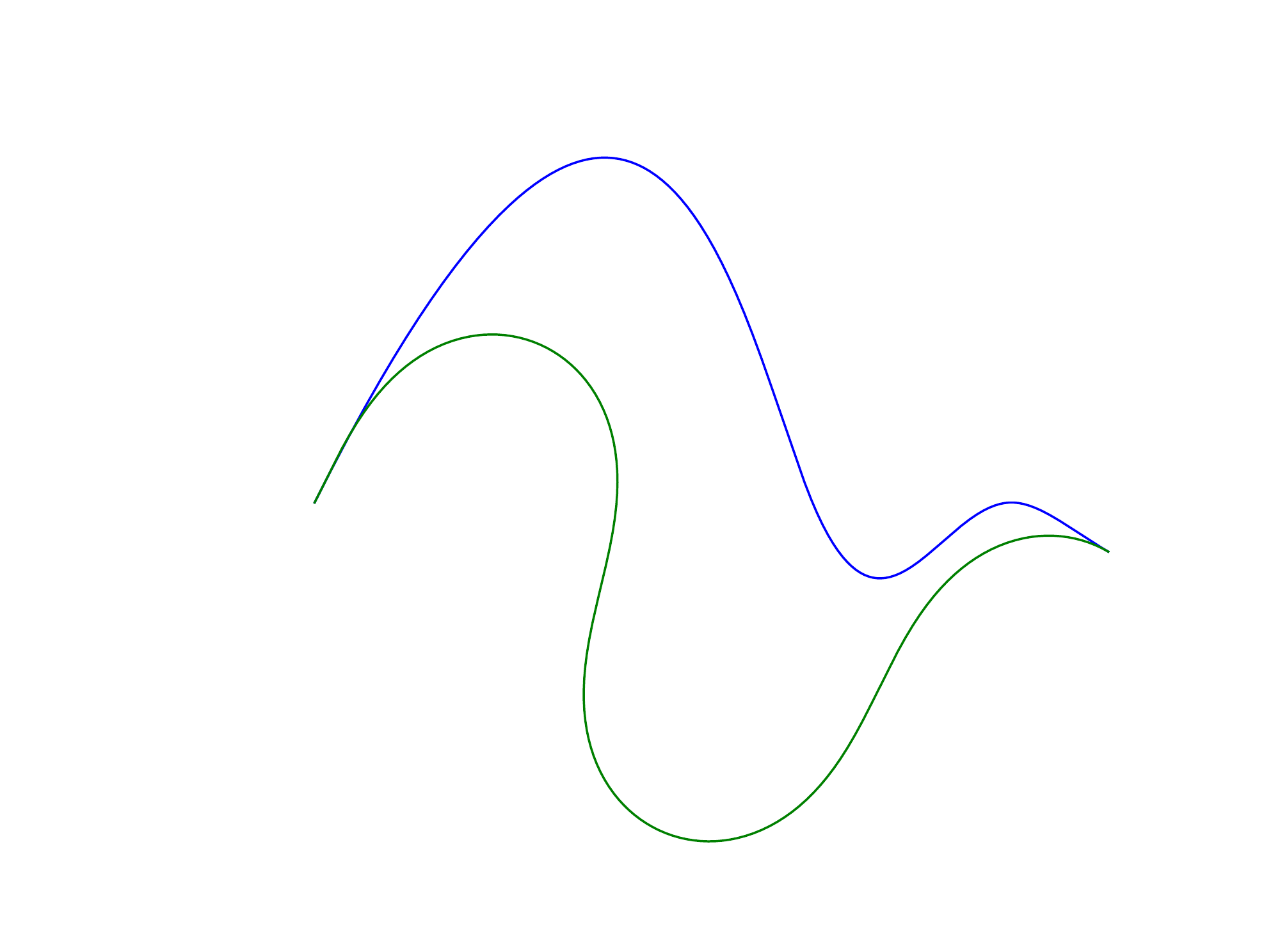}\hfil
\includegraphics[width=\unitlength,trim=40mm 20mm 20mm 35mm,clip]{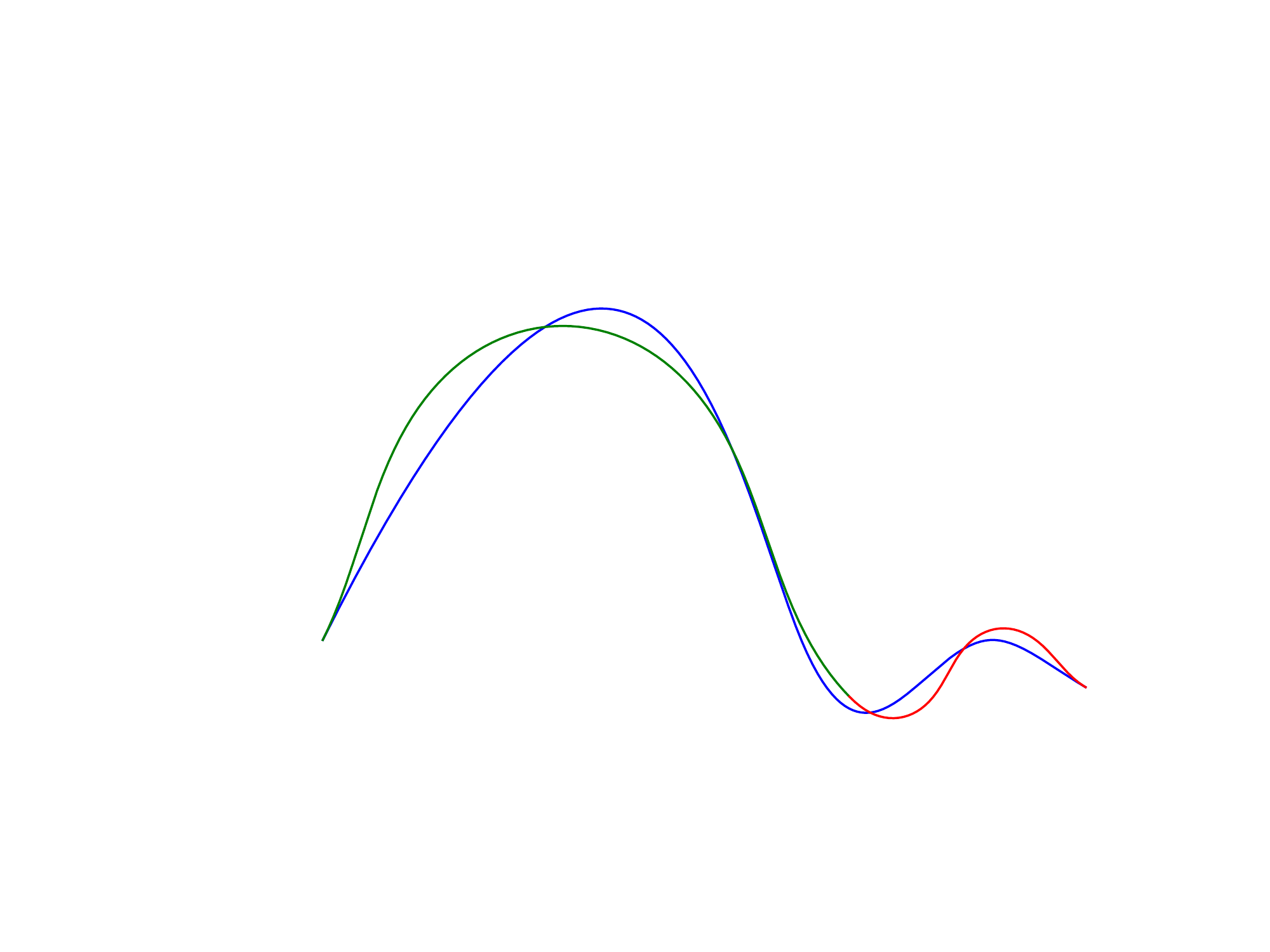}\hfil
\includegraphics[width=\unitlength,trim=40mm 20mm 20mm 45mm,clip]{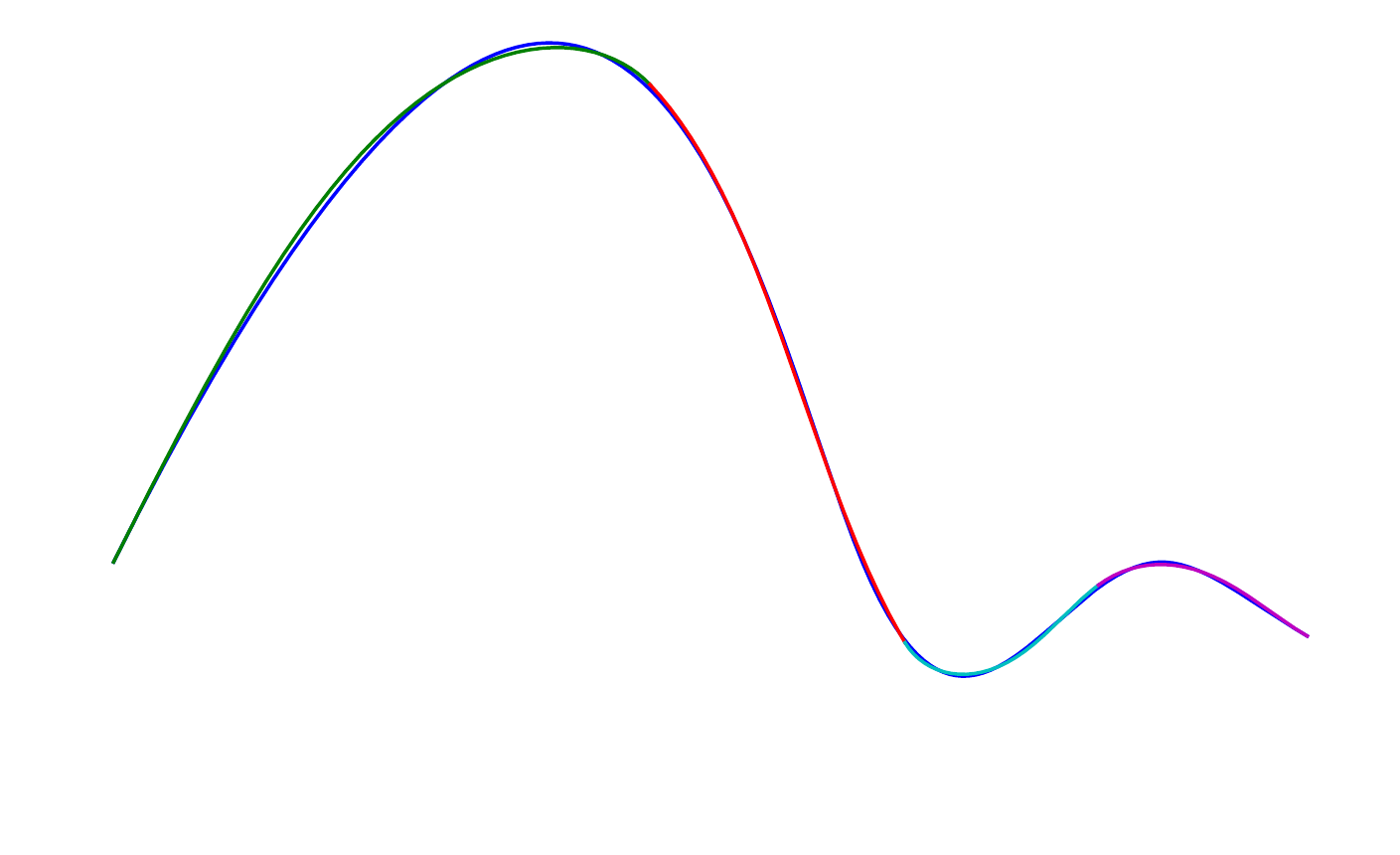}
\caption{Approximations of a more complex curve (blue) by, in order,
one, two and four tangent continuous  elastica segments. Both the endpoints and endtangents of
the target curve are matched.}
\label{fig:complex}
\end{figure}

It should be pointed out that the method for obtaining the initial guess (Section \ref{sec:elastic-parameters}) is only
useful for a curve segment that is not too far from some elastic curve segment.  Therefore, the practical use of this
 algorithm requires that a curve first be segmented into suitable pieces.  There are several ways to approach this:
for example (a) apply the initial guess method to the curve, (b)  measure the distance between the resulting elastic
segment and the original curve, (c) if the distance is too large, divide the curve into two and repeat.
An example with various segmentations is given in Figure \ref{fig:complex}.

In physical applications such as hot-blade cutting,
 it will also be necessary to consider the \emph{stability} of the elastic curve
segment obtained from this method; that is, small perturbations of the length, endpoints and tangents
of the curve should correspond to small changes in the solution shape.  See, for example,
\cite{lk} for a study of stability.  Constraints such as demanding that the curve segments have
no inflection points, or an upper bound on the curvature can be added to the procedure
to ensure stability.


\appendix
\section{Elliptic functions}\label{sec:app_elliptic}
We list some details of the Jacobi elliptic functions for convenience and to fix conventions.
Let $k\in (0,1)$. The \emph{elliptic functions} $\sn$, $\cn$ and $\dn$ with \emph{(elliptic) modulus} $k$ are defined as the solutions to the system of differential equations:
\begin{align*}\sn^\prime (u) &=\cn (u) \dn (u), & \sn (0)&=0, \\
\cn ^\prime (u) &=-\sn (u)\dn (u), & \cn (0)&=1, \\
\dn ^\prime (u) &=-k^2 \sn (u) \cn (u), & \dn (0)&=1.
\end{align*}
The \emph{complementary modulus} $k'\in[0,1]$ is defined by $k^2+k'^2=1$.
We have the identities:
\begin{equation}\label{eq:elliptic_identities}
 \sn^2u+\cn^2u=1, \quad \quad
  \dn^2u+k^2\sn^2u=1, \quad \quad
 \dn^2u-k^2\cn^2u=k'^2.
\end{equation}
The elliptic functions can be expressed in terms of integrals of trigonometric functions as follows.
Define the \emph{(elliptic) amplitude} $\am$ as:
\[
\am(t) = F^{-1}(t), \quad \quad F(\phi)=\int_0^\phi\frac1{\sqrt{1-k^2\sin^2u}} \du.
\]
 Then
 \begin{align*} \sn(u)=\sin(\am u), \quad \cn(u)=\cos(\am u), \quad \dn(u)=\sqrt{1-k^2\sin^2(\am u)}.\end{align*}

\subsection*{Elliptic integrals}
The integral $F(\phi,k)$ given, for each $k$, by the formula $F(\phi)$ above, is called
\emph{the incomplete elliptic integral of the first kind}.
 We define \emph{the incomplete elliptic integral of the second kind} by
\begin{align*}
 E(\phi,k)=\int_0^\phi\dn^2(u,k)\, \du.
\end{align*}
The \emph{complete} elliptic integrals of the first and second kind are respectively
the functions $F(k)=F(\frac\pi2,k)$ and $E(k)=E(\frac\pi2,k)$.

\subsection*{Addition formulas and periodicity}
 The elliptic functions satisfy the addition formulae:
\begin{align*}\sn(u+v)&=\frac{\sn u\cn v\dn v+\sn v\cn u\dn u}{1-k^2\sn^2u\sn^2v} \\ \cn(u+v)&=\frac{\cn u\cn v-\sn u\sn v\dn u\dn v}{1-k^2\sn^2u\sn^2v} \\ \dn(u+v)&=\frac{\dn u\dn v-k^2\sn u \sn v \cn u \cn v}{1-k^2\sn^2u\sn^2v}.\end{align*}

We define the \emph{quarter period} $K$ by:
\[
\am(K)=\frac\pi2, \quad \hbox{i.e.} \quad K=F\left(\frac{\pi}{2}\right),
\]
so that $\sn K=1$, $\cn K=0$ and  $\dn K=k'$.
 Then one obtains the periodicity:
\[
	\sn(u+2K)=-\sn u, \quad \quad
\cn(u+2K)=-\cn u, \quad \quad\dn(u+2K)=\dn u.
\]

\subsection*{Extension of $k$-domain}
The elliptic functions, as we have defined them, are only valid for $k\in[0,1]$. However, by analytic continuation (see e.g. Lawden \cite{lawden}), the domain of $k$ may be extended. For $k>1$, the following identities hold for all $u\in\R$: \begin{align*} \sn(u,k)&=\tfrac1k\sn(ku,\tfrac1k), \\ \cn(u,k)&=\dn(ku,\tfrac1k), \\ \dn(u,k)&=\cn(ku,\tfrac1k), \\
E(u,k)&=kE(ku,\tfrac1k)+u(1-k^2).
\end{align*}
Observe that $K(k)\to\infty$ as $k\to1$, so we cannot extend $K$ continuously. We choose the extension
\[
K(k)=\tfrac1{2k}K\left(\tfrac1k \right), \quad k>1,
\]
which ensures that the period of $\cn$ is always $4K$. We stress that this is \emph{not} the analytic continuation of $K$, which in fact takes non-real values for $k>1$.

\section{Derivatives}\label{sec:app_derivatives}
In this section, we list the derivatives of the basic elastica $\infl$, using the shorthand notation
\[
S=\sn(s,k),\ C=\cn(s,k),\ D=\sn(s,k),\ E=E(s,k).
\]
We have
\[
\infl(s,k)=\begin{pmatrix}2E-s \\ 2k(1-C)\end{pmatrix}.
\]

The derivatives with respect to $s$ follow from the definitions of the elliptic functions.
\begin{align*}
    \frac{\partial}{\partial s}\infl(s,k)&=\begin{pmatrix}2D^2-1 \\ 2kSD\end{pmatrix}, \\
    \frac{\partial^2}{\partial s^2}\infl(s,k)&=2kC\begin{pmatrix}-2kSD \\ D^2-k^2S^2\end{pmatrix} =2kC\begin{pmatrix}-2kSD \\ 2D^2-1\end{pmatrix}.
\end{align*}
The derivatives with respect to $k$ of $\sn$, $\cn$, $\dn$ and $E$ can be found in \cite{lawden}. From these,
with repeated use of \eqref{eq:elliptic_identities}, one can find
\begin{align*}
    \frac{\partial}{\partial k}\infl(s,k)=&\frac2{k'^2}\begin{pmatrix}k\left(SCD-EC^2-sk'^2S^2\right) \\ k'^2+C(k^2-D^2)-SD(E-sk'^2) \end{pmatrix}, \\
   \frac{\partial^2}{\partial s\partial k}\infl(s,k)=&\frac2{k'^2}\left(SD-C(E-sk'^2)\right)\begin{pmatrix}-2kSD \\ 2D^2-1\end{pmatrix}, \\
  \frac{k'^4}{2}   \frac{\partial^2}{\partial k^2}\infl(s,k) =& 	\begin{pmatrix}
    2SDC\left(D^2-k^2E^2+k'^2\left(s^2k^2-(E-s)^2-\tfrac12\right) \right)\\
\tfrac1k\left((1-2k^2S^2)(E-s)(2sk^2+E-s)C+DS(sk'^2-E)(4k^2C^2+k'^2) \right) \end{pmatrix} \\
  & +  \begin{pmatrix}
  (E-s)(C^2+D^2-4C^2D^2)+2sk^2(2S^2-1)D^2-sk'^2 \\
 -sk k'^2DS+s^2k^3C+kCS^2(2-2s^2k^4-2k^2S^2+k^2) \end{pmatrix}
\end{align*}

The derivatives of $\y_{(k,\stp,\ell,\w,\phi,x_0,y_0)}$ with respect to the control parameters can be found by straightforward calculations using the above.



\begin{thebibliography}{10}
\providecommand{\url}[1]{{#1}}
\providecommand{\urlprefix}{URL }
\expandafter\ifx\csname urlstyle\endcsname\relax
  \providecommand{\doi}[1]{DOI~\discretionary{}{}{}#1}\else
  \providecommand{\doi}{DOI~\discretionary{}{}{}\begingroup
  \urlstyle{rm}\Url}\fi


\bibitem{birkhoffdeboor1965}
Birkhoff, G., Boor, C.D.: Piecewise polynomial interpolation and approximation.
\newblock Approximation of Functions, Proc. General Motors Symposium 1964,
  H.L.~ Garabedian, ed. pp. 164--190 (1965).
\newblock Elsevier, Publ. Co. Amsterdam.

\bibitem{borbelyjohnson2014}
Borb\'ely, A., Johnson, M.: Elastic splines {I}: {E}xistence.
\newblock Constr. Approx. \textbf{40}, 189--218 (2014).

\bibitem{madifa}
Brander, D., B\ae{}rentzen, A., Evgrafov, A., Gravesen, J., Markvorsen, S.,
  N\o{}rbjerg, T., N\o{}rtoft, P., Steenstrup, K.: Hot blade cuttings for the
  building industry.
\newblock Preprint.

\bibitem{euler}
Euler, L.: Methodus inveniendi lineas curvas maximi minimive proprietate
  gaudentes; {A}dditamentum {I}: de curvis elasticis  (1744).
\newblock Translation: \cite{euler_translated}.

\bibitem{golombJerome1982}
Golomb, M., Jerome, J.: Equilibria of the curvature functional and manifolds of
  nonlinear interpolating spline curves.
\newblock SIAM J. Math. Anal \textbf{13}, 421--458 (1982).

\bibitem{horn1983}
Horn, B.: The curve of least energy.
\newblock ACM Trans. Math. Software \textbf{9}, 441--460 (1983).

\bibitem{lawden}
Lawden, D.: Elliptic Functions and Applications, \emph{Applied Mathematical
  Sciences}, vol.~80.
\newblock Springer-Verlag, New York (1989).

\bibitem{Levienthesis}
Levien, R.: From spiral to spline; optimal techniques for interactive curve
  design.
\newblock Ph.D. thesis, UC Berkeley (2009).

\bibitem{lk}
Levyakov, S.~V., Kuznetsov, V.~V.: Stability analysis of planar equilibrium configurations of
elastic rods subjected to end loads.
\newblock Acta Mech. \textbf{211}, 73--87 (2010).

\bibitem{Love}
Love, A.: A treatise on the mathematical theory of elasticity.
\newblock Cambridge University Press (1906).

\bibitem{malcolm1977}
Malcolm, M.: On the computation of nonlinear spline functions.
\newblock SIAM Journal on Numerical Analysis \textbf{14}, 254--282 (1977).

\bibitem{mehlum1974}
Mehlum, E.: Nonlinear splines.
\newblock Computer aided geometric design (Proc. Conf., Univ. Utah, Salt Lake
  City, Utah, 1974) pp. 173--207 (1974).

\bibitem{mumford1990}
Mumford, D.: Elastica and computer vision.
\newblock In: Algebraic geometry and its applications (West Lafayette, IN,
  1990), pp. 491--506. Springer, New York (1994).

\bibitem{euler_translated}
Oldfather, W.A., Ellis, C.A., Brown, D.M.: Leonhard {E}uler's elastic curves.
\newblock Isis \textbf{20}(1), pp. 72--160 (1933).
\newblock \urlprefix\url{http://www.jstor.org/stable/224885}

\bibitem{saalschutz}
Saalsch\"utz, L.: Der belastete Stab unter Einwirkung einer seitlichen Kraft.
\newblock B. G. Teubner, Leipzig (1880)

\bibitem{robarch}
S{\o}ndergaard, A., Feringa, J., N{\o}bjerg, T., Steenstrup, K., Brander, D., Gravesen, J., Markvorsen, S.,
B{\ae}rentzen, A., Petkov, K., Hattel, J., Clausen, K., Jensen, K., Knudsen, L., Kortbek, J.:
Robotic hot-blade cutting. 
\newblock In:  Robotic Fabrication in Architecture, Art and Design 2016 (pp. 150-164). Springer International Publishing.

\bibitem{ipopt}
A.~W\"achter and L.T. Biegler.
\newblock On the implementation of an interior-point filter line-search
algorithm for large-scale nonlinear programming.
\newblock {\em Math. Program., Ser. A} \textbf{106}, pages 25--57, 2006.

\end{thebibliography}

\end{document}